\documentclass[11pt]{amsart}

\usepackage{enumerate,url,amssymb,  mathrsfs}%, pdfsync}% ,showkeys, pdfsync}
\usepackage{graphicx}
\newtheorem{theorem}{Theorem}[section]
\newtheorem{lemma}[theorem]{Lemma}
\newtheorem*{lemma*}{Lemma}
\newtheorem{proposition}[theorem]{Proposition}

\theoremstyle{definition}
\newtheorem{definition}[theorem]{Definition}
\newtheorem{example}[theorem]{Example}
\newtheorem{conjecture}[theorem]{Conjecture}

\theoremstyle{remark}
\newtheorem{remark}[theorem]{Remark}

\numberwithin{equation}{section}

\newcommand{\abs}[1]{\lvert#1\rvert}

\newcommand{\A}{\mathbb{A}}
\newcommand{\C}{\mathbb{C}}

\newcommand{\DD}{\mathbb{D}}

\newcommand{\LL}{\mathcal{L}}
\newcommand{\R}{\mathbb{R}}
\newcommand{\Z}{\mathbb{Z}}
\newcommand{\T}{\mathbb{T}}
\newcommand{\Ho}{\mathscr{H}}
\newcommand{\Ha}{\mathcal{H}}

\newcommand{\const}{\mathrm{const}}
\newcommand{\dtext}{\textnormal d}

\DeclareMathOperator{\re}{Re}
\DeclareMathOperator{\im}{Im}

\DeclareMathOperator{\Mod}{Mod}
\def\Xint#1{\mathchoice
{\XXint\displaystyle\textstyle{#1}}%
{\XXint\textstyle\scriptstyle{#1}}%
{\XXint\scriptstyle\scriptscriptstyle{#1}}%
{\XXint\scriptscriptstyle\scriptscriptstyle{#1}}%
\!\int}
\def\XXint#1#2#3{{\setbox0=\hbox{$#1{#2#3}{\int}$}
\vcenter{\hbox{$#2#3$}}\kern-.5\wd0}}

\def\dashint{\Xint-}
\def\le{\leqslant}
\def\ge{\geqslant}
\begin{document}

\title{The Nitsche Conjecture}

\author{Tadeusz Iwaniec}
\address{Department of Mathematics, Syracuse University, Syracuse,
NY 13244, USA and Department of Mathematics and Statistics,
University of Helsinki, Finland}
\email{tiwaniec@syr.edu}
\thanks{Iwaniec was supported by the NSF grant DMS-0800416 and the Academy of
Finland grant 1128331.}

\author{Leonid V. Kovalev}
\address{Department of Mathematics, Syracuse University, Syracuse,
NY 13244, USA}
\email{lvkovale@syr.edu}
\thanks{Kovalev was supported by the NSF grant DMS-0913474.}

\author{Jani Onninen}
\address{Department of Mathematics, Syracuse University, Syracuse,
NY 13244, USA}
\email{jkonnine@syr.edu}
\thanks{Onninen was supported by the NSF grant  DMS-0701059.}

%    General info

\date{November 1, 2009}


\begin{abstract}
The Nitsche conjecture is deeply rooted in the theory  of doubly connected minimal surfaces. However,  it is commonly formulated in slightly greater generality as a question of existence of a harmonic homeomorphism between circular annuli
\[h\,\colon \;\mathbb A = A(r,R) \;\;\; \overset{\textnormal{\tiny{onto}}}{\longrightarrow}  \;\;\;A(r_\ast, R_\ast) =\mathbb A^* \]
In the  early 1960s, while attempting to describe all doubly connected minimal graphs over a given annulus $\mathbb A^*$, J.C.C. Nitsche observed that their conformal modulus  cannot be too large. Then he conjectured, in terms of isothermal coordinates, even more;  \\ \\
\begin{tabular}{|p{9.7cm}|}
\hline \vskip0.1cm
\begin{center}\textsl{A harmonic homeomorphism $h\colon \A \overset{\textnormal{\tiny{onto}}}{\longrightarrow}  \A^\ast$  exists if an only if:}\end{center}\vskip0.2cm
$$\frac{R_\ast}{r_\ast} \ge \frac{1}{2} \left(\frac{R}{r}+ \frac{r}{R}\right)$$
 \\\hline
\end{tabular}
\vskip0.5cm

This fascinating and engaging problem remained open for almost a half of a century. In the present paper we provide, among further generalizations, an affirmative answer to his conjecture.
\end{abstract}
\maketitle
 \tableofcontents
\section{Introduction and Overview}

Minimal surfaces in the Euclidean space  $\mathbb R^3$ are fundamental forms in nature, mathematics and physics in particular~\cite{DHKW,DHKW2, Os1}. Their elegant shape is due to the local minimum area property which yields zero mean curvature everywhere. One special class of surfaces arises by experimenting with the soap bubbles framed between two or more Jordan curves~\cite{Ba, Kau, Ni3}. For example, the minimal surface joining two coaxial circles in parallel planes has the shape of a catenoid, a configuration that is extremal for numerous problems~\cite{Ba, Hi,  Kau, Ni1, Ni3, Ni4, OsSc}. Theorem~\ref{modS} in the present paper shows that catenoid has also the largest conformal modulus among minimal graphs over a given annulus.
\begin{definition}
 An open doubly-connected surface  $\mathbb  S $ in $\mathbb R^3$  is a conformal image of either
\begin{itemize}

\item{punctured complex plane  $\mathbb C_\circ  =  \{ z \in \mathbb C ;\; z \neq 0\}$}
\item {punctured disk $\mathbb D_\circ = \{ z \in\mathbb C ;\; 0< |z| < 1\}$}
\item {or an annulus  $\mathbb A = A(r,R) =   \{z;\; r<|z| < R \}\;$ ,  $\;\;\;0<r <R< \infty$ }

\end{itemize}
\end{definition}

In this latter case the \textit{conformal modulus} of  $\mathbb S$  is defined by
\begin{equation}
    \textrm {Mod}  \mathbb \;\mathbb S  =  \log \frac{R}{r}   > 0
\end{equation}

Note that the ratio $ \frac{R}{r} $ is independent of the choice of the conformal mapping, say  $F = (u,v,w) \colon \mathbb A  \overset{\textnormal{\tiny{onto}}}{\longrightarrow} \mathbb S $.  This is immediate from the classical result of Schottky~\cite{Sc},  1877.
Let us first establish the notation of two circular annuli in the complex plane that will remain standard throughout this paper

\begin{equation}\nonumber
\begin{split}
\begin{cases}\A \;\;= A(r,R)\;\;\;= \left\{ z\in \C \colon \;r\;< \abs{z}\;<R \;\, \right\}\, , \qquad 0<r\;<R\;\,< \infty \\
\A^\ast = A(r_\ast,R_\ast)= \left\{ \xi \in \C \colon r_\ast< \abs{\xi}<R_\ast  \right\}\, , \qquad 0<r_\ast<R_\ast< \infty
\end{cases} \end{split}
\end{equation} \\
A theorem  of Schottky asserts that:

 \begin{theorem}\label{Schottky}
  There exists a conformal homeomorphism $h \colon \A  \overset{\textnormal{\tiny{onto}}}{\longrightarrow} \A^\ast $  if and only if the annuli have the same modulus; that is, $$\Mod \A= \log \frac{R}{r}  =  \log \frac{R_\ast}{r_\ast }  =  \Mod \A^\ast$$
  Moreover, up to the rotation of the annuli, every such map takes the form\\
  \begin{equation}\nonumber
  \begin{split}
  h(z) =  \begin{cases} \frac{r_\ast z}{r}  \,, \;\;\;\;\;\textrm{ if preserving the order of boundary components} \\\;\\ \frac{ rR_\ast}{ \, z}\;, \;\;\;\;\;\textrm{if reversing the order of boundary components}
  \end{cases}\end{split}
  \end{equation}
\end{theorem}

We shall later devote Section~\ref{Sec3.1} to a fresh proof of Schottky's theorem.  It is this proof that reveals the novelty of our approach and, to some extent, leads to computational advances needed for harmonic mappings of annuli.

The coordinate functions $\,F = (u,v,w)  \colon \A  \overset{\textnormal{\tiny{onto}}}{\longrightarrow} \mathbb S \,$ in the conformal representation of a surface, called \textit{isothermal parameters}\,, satisfy the Cauchy-Riemann system:

 \begin{equation}\nonumber
 \begin{split}
  \begin{cases} u_xu_y  +v_xv_y + w_x w_y  =  0  \,, \;\; \;\;\;\;\;\;\;\textrm{-\;vectors} \;\; F_x \;\textrm{and} \; F_y \; \textrm{are orthogonal in} \;\;\mathbb R^3\\
  u_x^2  + v_x ^2 + w_x^2  =  u_y^2 + v_y^2 +w_y^2  \,,\;\;\;\;\; \;\textrm{-\;vectors}\;\;  F_x \;\textrm {and} \;F_y\;\;\textrm{have equal length ,}
 \end{cases}\end{split}
 \end{equation}
Throughout this paper we take advantage of the complex variable $\, z = x + i y\,$ and the two Wirtinger differential operators $$\frac{\partial}{\partial z} = \frac{1}{2} \left(\frac{\partial}{\partial x} - i \frac{\partial}{\partial y}\right )\; \;\;\;\textrm {and}\;\;\;\;  \frac{\partial}{\partial \bar z} = \frac{1}{2} \left(\frac{\partial}{\partial x} + i \frac{\partial}{\partial y}\right) \;, \;  $$
In this notation the Cauchy-Riemann system takes the form of one complex equation
 $$
 u_z^2  + v_z ^2 + w_z^2  = 0\,,\;\;\; \textrm{where}\;\;  (u_z, v_z, w_z) =   F_z = \frac{\partial}{\partial z} F
 $$
Now, the surface is minimal if and only if the isothermal parameters are harmonic or, equivalently, the complex vector field $ F_z \colon \mathbb A \rightarrow  \mathbb C^3 $ is holomorphic; that is, $$ \frac{\partial}{\partial \bar z} \left(F_z\right)  =  F_{z\bar z} =  \frac{1}{4} \Delta F  =  0 .$$
Let us factor the ambience of the surface into the complex plane and the real line $\mathbb S \subset \mathbb R^3 \simeq \mathbb C\times \mathbb R$. Thus $F = (h, w) \colon \mathbb A \rightarrow  \mathbb C\times \mathbb R $,  where $ h =  u+ iv \colon \mathbb A \rightarrow  \mathbb C $ is a complex harmonic map. A simple direct computation shows that
$$
  u_z^2  +  v_z ^2  \, = \, h_z \,\overline{h_{\bar z}}\;  =  - w_z^2
$$
Hence, all zeros of the holomorphic function $ \, h_z \overline{h_{\bar z}}\;$ must have even order. This allows us to determine the third isothermal parameter in terms of $h$
$$
   w  =  \im \int \sqrt{\, h_z \,\overline{h_{\bar z}}\;} \;\textrm{d} z
$$
It is evident that every complex harmonic map $ h =  u+ i\,v \colon \mathbb A \rightarrow  \mathbb C $, for which $ \, h_z \overline{h_{\bar z}}\;$ admits a continuous branch of square root, can be lifted to the isothermal parameters of a minimal surface. The surface is flat; that is,  $ w \equiv const $, if and only if $h$ is holomorphic (orientation preserving) or antiholomorphic (orientation reversing).

Let us consider now the case in which  $\,\mathbb S\,$ is a minimal graph over an annulus $\mathbb A^* = A(r_*, R_* ) \subset \mathbb C$. Thus
$h \colon A(r,R)  \overset{\textnormal{\tiny{onto}}}{\longrightarrow} \mathbb A^{\ast}\;$ is a harmonic homeomorphism. Without losing any generality we can assume that $h$ preserves the orientation. This yields an elliptic Beltrami type equation for $h :$
\begin{equation}
\overline{h_{\bar z}}  =  \mu(z)\, h_z \;,\;\;\textrm{with} \;\;  |\mu(z)| < 1\;,\;\;\textrm{and} \;\;\;w_z  =  \pm\;i\,\sqrt{\mu} \; h_z\,,
\end{equation}
where $\,\sqrt{\mu}\,$ exists as a single valued analytic branch in $\mathbb A = A(1,R)$. Our basic example of this is an upper slab of a catenoid furnished by
the parameters:
\begin{equation}
  h(z) =  \frac{1}{2} \left( z + \frac{1}{\bar z}\right ) ,\;\;\mu (z) = \frac{- 1}{z^2}\;,\;\;\;  w(z)  = \log |z| \,,\;\;\; 1 <|z| < R.
\end{equation}
This is a minimal graph over the annulus $\mathbb A^*= A(1,R_*)\,,$  $\; R_* =  \frac{1}{2} (R + \frac{1}{R})$
\begin{equation}
\mathbb S^* =  \left\{ (\xi,w) \in \mathbb A^*\times \mathbb R\,; \;\; \; w = \log \,(\, |\xi|  +  \sqrt{|\xi|^2 \;- 1} \,)\;<\log R \right\}
\end{equation}
 Note that the width of the slab $\mathbb S^*$, distance from the base circle to the top one,  equals its conformal modulus, $$\Mod \;\mathbb S^*  =  \log R  = \log\,(\,R_* +  \sqrt{R_*^2 - 1} ) .$$ Suppose we are now given two frames (Jordan curves) to be used to create a soap bubble of a doubly connected minimal surface. Generally speaking, when these frames are moved farther apart, the conformal modulus of a surface increases and there is a moment when the bubble brakes down. The critical upper bound of the conformal modulus for which the minimal graph exists is the essence of the Nitsche conjecture, which is now one of our results.
\begin{theorem}\label{modS}
Let $\mathbb S$ be a minimal graph over the annulus $\A^\ast=A(r_\ast, R_\ast)$, then
\begin{equation}
\Mod \; \mathbb S \; \leqslant \;\log \left (\frac{R_\ast}{r_\ast}  +  \sqrt{\frac{R_\ast^2}{r_\ast^2} \;-\;1} \right)
\end{equation}
Equality occurs if and only if $\;\mathbb S\;$ is the upper slab of a catenoid.
\end{theorem}

\begin{center}
\begin{figure}[h]
\includegraphics[height=70mm]{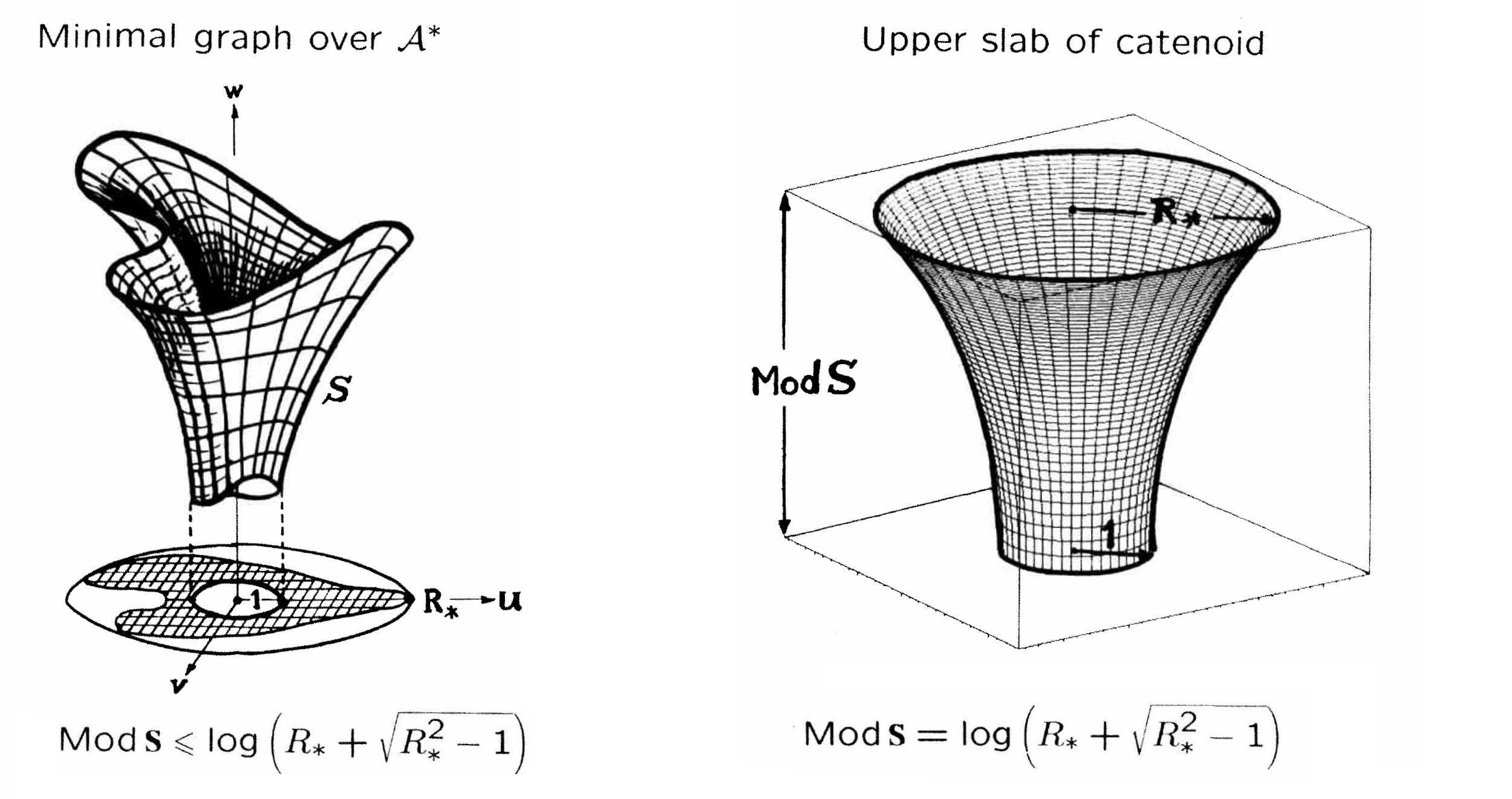}
\caption{Among all minimal graphs over given annulus the upper slab of catenoid has the greatest conformal modulus.}
\end{figure}
\end{center}

In 1962 J.C.C. Nitsche announced, in a short article ~\cite{Ni}, that the existence of a harmonic homeomorphism
$h \colon \A  \overset{\textnormal{\tiny{onto}}}{\longrightarrow} \A^\ast$, whether or not it comes from a minimal graph, yields a lower bound on $\Mod \A^\ast$ in terms of $\Mod \A$. He conjectured that the necessary and sufficient condition for  such a mapping to exist is  the following inequality, now known as the {\it Nitsche bound}
\begin{equation}\label{nb}
\frac{R_\ast}{r_\ast} \ge \frac{1}{2} \left(\frac{R}{r}+ \frac{r}{R}\right)
\end{equation}
Subsequently, his conjecture appeared in  monographs~\cite[\S 878]{Nib}, \cite[p.~138]{Dub}, \cite[Conj. 21.3.2]{AIMb} and surveys~\cite{BH,Lyz,Sch}. Various lower bounds for $R_\ast/r_\ast$ have been obtained by Lyzzaik~\cite{Ly}, Weitsman~\cite{We}, Kalaj~\cite{Ka}, and by Nitsche himself (see~\cite{Lyz}). Worth noting is Lyzzaik's bound in ~\cite{Ly} that exhibits linear growth of $\,\frac{R_{\ast}}{r_{\ast}}\,,$ as the ratio $\,\frac{R}{r}\,$ approaches infinity. However, none of these results came close to ~\eqref{nb}. Here we prove this inequality, in somewhat greater generality, which obviously implies Theorem ~\ref{modS}.

One should note in advance that whenever $h \colon \A \overset{\textnormal{\tiny{onto}}}{\longrightarrow} \A^\ast$ is a homeomorphism, the function $z \mapsto \abs{h(z)}$ extends continuously up to the boundary of~$\A$. Moreover, there are two possibilities; either
\begin{equation}\label{boundcond}
\begin{split}
\abs{h(z)}=  \begin{cases}r_\ast & \qquad \mbox{for} \quad \abs{z}=r \\
R_\ast & \qquad \mbox{for} \quad \abs{z}=R
\end{cases}\end{split}
\end{equation}
or  the other way round
\begin{equation}
\begin{split}
\abs{h(z)}=  \begin{cases}R_\ast & \qquad \mbox{for} \quad \abs{z}=r \\
r_\ast & \qquad \mbox{for} \quad \abs{z}=R
\end{cases}\end{split}
\end{equation}
In the former case, we say that $h$ is consistent with the order of boundary components. We also have two additional possibilities depending on whether $h$ preserves or reverses the orientation. Accordingly, there are four homotopy classes of homeomorphisms between annuli. Without loss of generality we shall confine ourselves to studying the class
\begin{equation}
\begin{split}
\;\;\;\;\Ho (\A, \A^\ast) = \big\{ h \colon  \A  \overset{\textnormal{\tiny{onto}}}{\longrightarrow}  & \A^\ast \colon \mbox{sense-preserving homeomorphisms} \\ & \mbox{satisfying the boundary condition~\eqref{boundcond}} \big\}
\end{split}
\end{equation}
The following notation will represent harmonic mappings in this class,
\begin{equation}
\Ha (\A, \A^\ast) = \left\{h \in \Ho (\A, \A^\ast) \colon \Delta h=0   \right\}
\end{equation}

Our preeminent result, though not the most general, is:

\begin{theorem}\label{th1}\textnormal{\textsc{(Existence of Harmonic Maps)}}
The annuli $\A=A(r,R)$ and $\A^\ast=A(r_\ast, R_\ast)$ admit a harmonic homeomorphism $h \colon \A \overset{\textnormal{\tiny{onto}}}{\longrightarrow}  \A^\ast$ if and only if
\begin{equation*}
\frac{R_\ast}{r_\ast} \ge \frac{1}{2} \left(\frac{R}{r}+ \frac{r}{R}\right)
\end{equation*}
When equality occurs, then every $\,h\,\in \mathcal H(\mathbb A, \mathbb A^*)\,$ takes the form
\begin{equation*}
h(z)\; = \; \frac{1}{2} \left(\frac{z}{r}+ \frac{r}{\bar z}\right)r_\ast \,  e^{i \alpha}\, , \quad 0 \le \alpha < 2 \pi
\end{equation*}
\end{theorem}
Various geometric and analytic properties of harmonic mappings  are discussed, e.g.,  in the books ~\cite{Dub, EFb, Job} and articles~\cite{BH, BH1, HJW, HKW, Sch}.
There are notable recent studies~\cite{AIM, AIMb, IO} that also concluded with the Nitsche bound in~\eqref{nb}, though from a slightly different framework. This is the framework of the Dirichlet energy
\begin{equation}
\mathcal E [h]= \iint_\A \abs{Dh}^2 = 2 \iint_\A \left(\abs{h_z}^2+ \abs{h_{\bar z}}^2\right)
\end{equation}
In general, minimizing the energy among homeomorphisms need not lead to the Laplace equation. The reason is that passing to the weak limit of a minimizing sequence of homeomorphisms in $\Ho (\A, \A^\ast) $  may result in a noninjective mapping. Harmonicity is lost exactly at  the branch points near which the extremal mapping fails to be injective. Outside the branch set the extremal mappings are indeed harmonic. This latter fact follows from Rad\'o-Kneser-Choquet Theorem \cite[p. 29]{Dub}. The above recent studies of the minima of the Dirichlet energy can be encapsulated as:
\begin{theorem} %\label{th3aaa}
The class $\Ho (\A, \A^\ast)$ of homeomorphisms $h \colon \A \overset{\textnormal{\tiny{onto}}}{\longrightarrow}  \A^\ast$ admits an energy minimizer if and only if
\begin{equation*}%\label{nitschebound}
\frac{R_\ast}{r_\ast} \ge \frac{1}{2} \left(\frac{R}{r}+ \frac{r}{R}\right)
\end{equation*}
All minimizers are harmonic, in symbols:
\begin{equation}
  \inf_{h\in\mathscr H(\mathbb A\,,\,\mathbb A^*)} \iint_\A \abs{Dh}^2\;=\; \inf_{h\in\mathcal H(\mathbb A\,,\,\mathbb A^*)} \iint_\A \abs{Dh}^2\;=\; \iint_\A \abs{Dh^{\min}}^2\;
\end{equation}
Moreover, modulo rotation,  they take the form \begin{equation} h^{\min}(z)=az+b\bar z^{-1}, \;\; a = \frac{R R_* - r r_*}{R^2 -r^2}\,,\;\; b= \frac{(Rr^* -rR^*) rR}{R^2-r^2}\end{equation}
\end{theorem}
Note that the necessary condition is immediate from Theorem~\ref{th1}, though it was first established in ~\cite{AIM, IO} by completely different idea (the method of free Lagrangians). Indeed, beyond the Nitsche bound there is no homeomorphism minimizing the energy because otherwise it would be a harmonic map, by  Rad\'o-Kneser-Choquet Theorem \cite[p. 29]{Dub}.

It is interesting to know what happens beyond the inner boundary of $\, A(r,R)\,$. For simplicity, suppose that $\,r=r_{\ast} = 1\,,$ so in  the critical configuration we have $ \,R_*  = \frac{1}{2} (R + R^{-1})$. The Nitsche map $h(z) =  \frac{1}{2} (z + \bar{z}^{-1})$, so extended to the annulus $\, A(R^{-1}, R ) $, becomes a harmonic double cover of $A(1,R_*)$. Folding takes place along the unit circle where the Jacobian vanishes. On the other hand, passing to a weak limit of a  minimizing sequence of homeomorphisms $h \colon A(R^{-1}, R )\overset{\textnormal{\tiny{onto}}}{\longrightarrow}  \;\mathbb A^*$    results in squeezing the inner portion $\,A(R^{-1}, 1)\,$ of the annulus $\,A(R^{-1},\, R)\,$ onto the unit circle. The minimizer, found in ~\cite{AIM, IO}, is unique modulo rotation and takes the form:

\begin{equation}\nonumber
  \begin{split}
  h^{\inf}(z) =  \begin{cases} \frac{ z}{|z|}  \,, \;\;\;\textrm{if} \;\;R^{-1} < |z|<1 \;,\;\textrm{ -hammering phenomenon } \\\;\\ \frac{1}{2} (z+\bar z^{-1})\;, \;\;\textrm{if}\;\; |z|>1 \;,\;\;\textrm{-the critical Nitsche map}
  \end{cases}\end{split}
  \end{equation}

  The critical Nitsche map $h(z) =  \frac{1}{2} (z + \bar{z}^{-1})$, extended harmonically beyond the unit circle and  furnished with the vertical isothermal coordinate $ w =\log |z| \,, $ for $\;\frac{1}{R} <|z| < R\, $, gives rise to a symmetric slab of the catenoid. This result strongly suggests nonexistence of any harmonic homeomorphism $h \colon A(R^{-1}, R )\overset{\textnormal{\tiny{onto}}}{\longrightarrow}  \;\mathbb A^*$. But
  \emph{\`a priori}, there might exist harmonic homeomorphisms which do not minimize the energy. To show that this is not the case we found that the integral means
\begin{equation} \label{means}
  U(\rho) \;= \dashint_{\T_\rho} \abs{h}^2 = \frac{1}{2 \pi \rho} \int_{\abs{z}=\rho} \abs{h(z)}^2 \, \abs{\dtext z}
\end{equation}
over the circles
\begin{equation*}
\T_\rho= \left\{z \in \C \colon \abs{z}=\rho\right\}\, , \quad r \le \rho <R\,
\end{equation*}
are better suited than the energy integrals.

It should be observed that  harmonicity of a function $h=h(z)$ is invariant under conformal change of the $z$-variable. Therefore, the Nitsche bound remains valid for harmonic homeomorphisms  defined on any doubly connected domain whose conformal  modulus coincides with that of $\A$.  It is therefore of interest to look at the role of  the boundary curves in the target annulus as well. Although we  need the inner boundary of the target to be a circle, the circular shape of the outer boundary turns out to be  inessential, there remains a substitute of the Nitsche bound in terms of the integral means $U(\rho)$. In our generalized form of the Nitsche bound the target will be a half circular annulus; that is,  a doubly connected domain, denoted by  $\mathcal A^\ast$,  whose inner boundary is a circle $\T_\ast= \left\{w \in \C \colon \abs{w}=r_\ast\right\}$. We do not specify the outer boundary of $\mathcal A^\ast$ as it can be arbitrary. Let $\Ha (\A, \mathcal A^\ast)$ denote the class of orientation preserving harmonic homeomorphisms $h \colon \A  \overset{\textnormal{\tiny{onto}}}{\longrightarrow}  \mathcal A^\ast  $ which preserve the order of the boundary components. In particular, $\;\lim\limits_{\abs{z} \searrow \,r}\abs{h(z)}=r_\ast$.

\begin{center}
\begin{figure}[h]
\includegraphics[height=60mm]{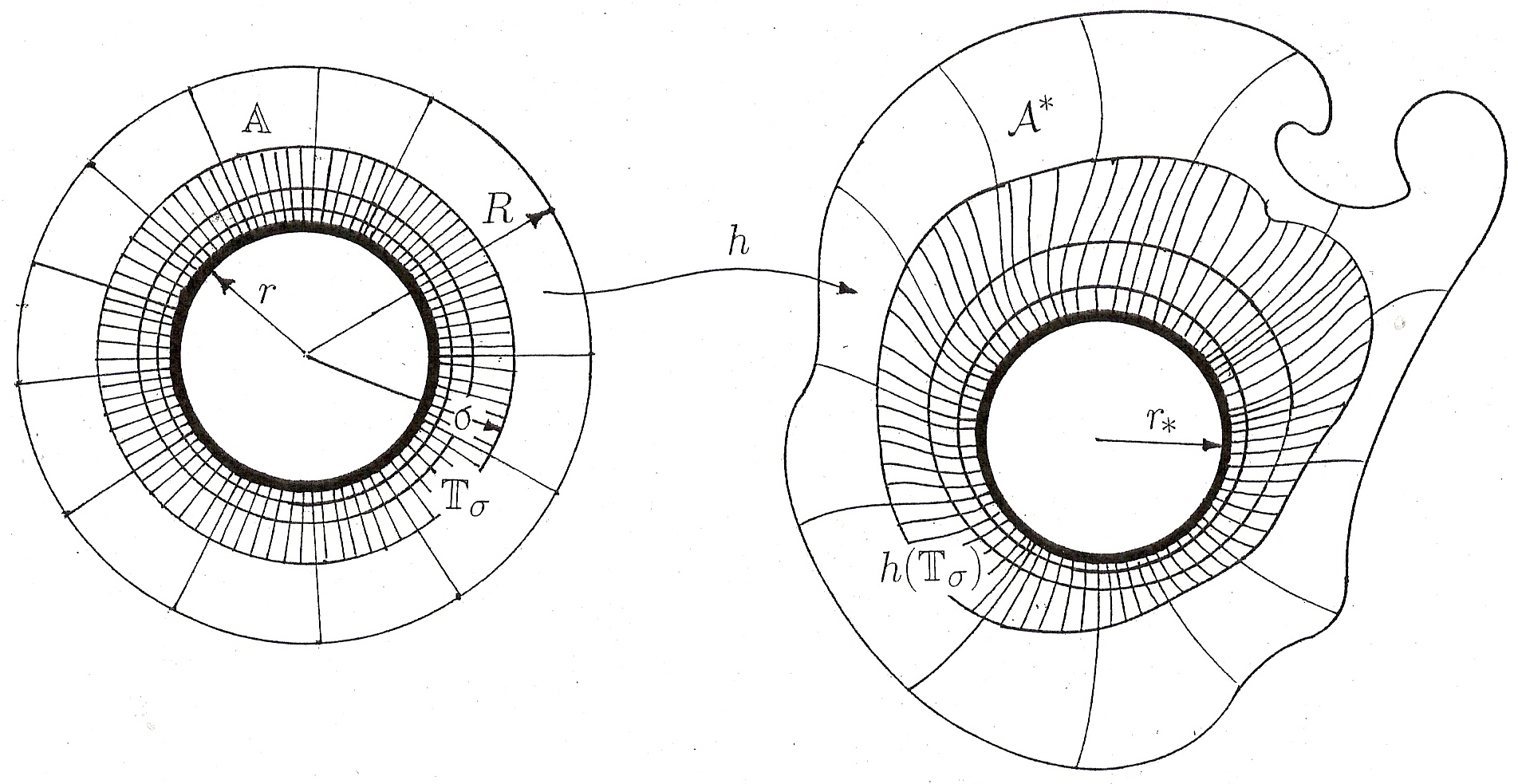}
\caption{Harmonic Evolution of Circles}
\end{figure}
\end{center}

\begin{theorem}\label{th2} \textnormal{\textsc{(Generalized Nitsche Bound)}} For every $h\in \Ha (\A, \mathcal A^\ast)\;$ we have
\begin{equation}\label{genbound}
\left[ \dashint_{\T_\sigma} \abs{h}^2  \right]^\frac{1}{2} \ge \frac{1}{2} \left(\frac{\sigma}{r}+ \frac{r}{\sigma}\right)r_\ast\; , \quad r \le \sigma <R
\end{equation}
If equality occurs at some radius $\sigma \in (r,R)$, then it holds for every $\sigma \in [r,R)$ and $h(z)= \frac{1}{2}\left(\frac{z}{r}+\frac{r}{\bar z}\right)r_\ast\,  e^{i \alpha}\,,$ for some $0\leqslant \alpha < 2\pi$.
\end{theorem}

Such a more general statement   not only strengthens  the Nitsche Conjecture, but also is the key to the proof. Theorem~\ref{th2} should be viewed as a sharp lower estimate for the growth of integral means of harmonic mappings under certain initial constraints. These constraints concern topological behavior of $h$ near the inner boundary $\T_r$ rather than its boundary values. It causes no loss of generality to assume that $h$ is smooth up to the inner boundary; we shall reduce ourselves to smooth mappings via Lemma~\ref{lemappr}. What one really needs for smooth mapping  $\,h\,$ is its harmonicity in $\A(r,R)\,$ and the following three initial conditions:\\

\begin{enumerate}[(I)]
\item\label{in1} $\;\;\;\;\;h \colon \T_r \to \T_r\;\;\;$ is a homeomorphism homotopic to the identity;
\item\label{in2} $\;\;\;\;\;\displaystyle \frac{\dtext}{ \dtext \sigma} \dashint_{\T_\sigma} \abs{h}^2 \ge 0\;,\; $ at $\;\sigma =r$;
\item\label{in3} $\;\;\;\;\;\displaystyle \int_{\T_r} \det Dh \ge 0 $.
\end{enumerate}
In fact these three conditions yield inequality~\eqref{genbound}. Even more, it turns out that condition~\eqref{in3} is redundant when  $\Mod \A \le 1$. The proof of this case of Theorem~\ref{th2} is based upon the ideas from our earlier paper~\cite{IKO}. However, Example~\ref{ex} demonstrates that constraints~\eqref{in1}--\eqref{in2} are insufficient to deduce inequality~\eqref{genbound} when $\Mod \A$ is large. A new approach is required in which we make use of the Jacobian bound~\eqref{in3}. The key ingredient is the following result concerning harmonic self-mappings of the unit disk, which is interesting in its own right.
\begin{theorem}\label{th3}
Let a harmonic homeomorphism $f \colon \overline{\mathbb D}  \overset{\textnormal{\tiny{onto}}}{\longrightarrow} \overline{\mathbb D}$ be $\mathscr C^1$-smooth in the closed unit disk $\,\overline{\mathbb D} = \{z \in \C \colon \abs{z} \le 1\}$. Then
\begin{equation}\label{th3form}
\int_{\partial \mathbb D} \abs{\det D\!f} \ge \iint_{\mathbb D} \abs{D\!f}^2 \ge 2 \iint_{ \mathbb D} \abs{\det D\!f}= 2 \pi
\end{equation}
The first inequality is strict unless  $f$ is an isometry.
\end{theorem}

We end this overview by stating an analogue of Nitsche's conjecture for noninjective harmonic maps, which has geometric interpretation for more general doubly connected minimal surfaces, not necessarily graphs.

\begin{conjecture} Let $\;h\,\colon A(r,R) \, \overset{\textnormal{\tiny{into}}}{\longrightarrow}  \,A(r_\ast, R_\ast) $
be a harmonic map (not necessarily injective) with nonzero winding number; that is,
$$
\int_{\T_\sigma} \frac{d h}{h} \ne 0\;, \;\;\;\textrm{for some (equivalently,  for all)}\;\;\; \sigma \in (r,R)
$$
Then
\begin{equation}
\frac{R_\ast}{r_\ast} \ge \frac{1}{2} \left(\sqrt{\frac{R}{r}}+ \sqrt{\frac{r}{R}}\right)
\end{equation}
Equality occurs for the double cover Nitsche map:
\begin{equation}
    h(z)  =  \frac{1}{2} \left( \frac{z}{\sqrt{rR}}  \;+\;\frac {\sqrt{rR}}{\bar z}\right)\;,\;\;\;\; r<|z|<R
\end{equation}
\end{conjecture}
The solution of the Nitsche conjecture and the methods presented here open new questions yet to be answered.

\section{Preliminaries}
The annuli $\A$ and $\,\mathbb \A^\ast$ as well as the half circular annulus $\mathcal A^\ast$ will henceforth be rescaled so their inner boundaries are the unit circles, still denoted by $\mathbb T.$  Whenever the target of $h \in \Ha (\A, \mathcal A^\ast)$ need not be specified we abbreviate this notation to $h \in \Ha (\A, \ast)$. Thus, from now on
$$
r=r_\ast=1\;, \;\;\;\;\;\; |h(z)| \equiv 1 \;\;\;\text{for}\;\; z\in \mathbb T
$$

The complex derivatives $h_z$ and $h_{\bar z}$ in polar coordinates are:
\begin{equation*}
h_z = \frac{1}{2}  \left(h_\rho - \frac{i}{\rho}h_\theta  \right)e^{-i\theta}\, , \quad h_{\bar z} = \frac{1}{2} \left(h_\rho + \frac{i}{\rho}h_\theta  \right)e^{i\theta} , \quad z=\rho e^{i \theta}
\end{equation*}
Hence one finds the Hilbert-Schmidt norm of the differential matrix
$$
 |Dh|^2  = 2 ( |h_z|^2 + |h_{\bar z}|^2 )  =   |h_\rho|^2  +  \rho^{-1} |h_\theta|^2 \;,$$
the Laplacian
$$\Delta h  =  4 \frac{\partial^2 h}{\partial z\partial \bar z} =  \frac{1}{\rho} \frac{\partial}{\partial\rho} \left(\rho\frac{\partial h}{\partial \rho}\right)  +\;\frac{1}{\rho^2}\frac{\partial^2 h}{\partial \theta^2}\;=\; h_{\rho\rho} \,+\, \frac{1}{\rho} h_\rho \,+\, \frac{1}{\rho^2} h_{\theta\theta}\,
$$
and the Jacobian determinant
$$
  J(z,h) \,=\, |h_z|^2\,-\, |h_{\bar z}|^2  =  \frac{\im ({\bar h_\rho}\, h_\theta )}{\rho} \,.
$$
We shall work with a number of circular means to which the following commutation rule will apply
\begin{equation}\label{comrule}
\frac{\dtext}{\dtext \rho} \dashint =  \dashint \frac{\dtext}{\dtext \rho} \, , \qquad \qquad 1 < \rho < R
\end{equation}
The main objects are the circular means of $\,|h|^2\,$ for mappings  $h \in \Ha (\A, \ast)$,
\begin{equation*}
U(\rho)= \dashint_{\T_\rho} \abs{h}^2 \, , \qquad U(1) = 1 \,, \qquad \dot{U}(\rho)= \dashint_{\T_\rho} \abs{h^2}_\rho\, , \qquad \qquad 1 <  \rho < R
\end{equation*}
We shall also discuss so-called mean  radius of $\,h(\mathbb T_\rho)\,,$ defined by the rule:
\begin{equation*}
\mathscr R_h(\rho) \,=\, \mathscr R(\rho)=  \sqrt{U(\rho)} \;=\;\left(\dashint_{\T_\rho} \abs{h}^2\right)^{\frac{1}{2}} \, ,\;\;\mathscr R(1) = 1 \,,\;\; \dot{\mathscr R}=  \;\dot U\,/2\,\sqrt{U}
\end{equation*}
Here, as usual, dot over $U$ and $\,\mathscr R\,$ stands for the $\rho$-derivative.
Thus, we aim to prove the inequality
\begin{equation*}
\mathscr R(\rho) \ge  \frac{1}{2}\left( \rho  + \frac{1}{\rho}\right) \, , \qquad 1 \le \rho <R
\end{equation*}
 Next observe that away from the outer boundary of $\mathbb A$ the energy of $h \in \Ha (\A, \ast)$ is always finite.
The proof is an exercise with Green's formula; first we have for $1<\varepsilon<\rho <R\,,$
\begin{equation}\label{eq25}
\begin{split}
\frac{1}{\pi} \iint_{A(\varepsilon,\, \rho)} \abs{Dh}^2 \;&=\;\frac{1}{2\pi} \iint_{A(\varepsilon, \,\rho)} \Delta\abs{h}^2 {}  \\
 &= \frac{1}{2\pi} \int_{\partial\, A(\varepsilon, \,\rho)} \abs{h^2}_{\textsc{n}} \; =\;\rho\,\dot{U}(\rho) \; -\; \varepsilon\,\dot{U}(\varepsilon) {}
\end{split}
\end{equation}
Hence we infer that $\rho\,\dot{U}(\rho)$ is strictly increasing and then
\begin{equation*}
\begin{split}
 ( \rho \log \rho ) \;\dot{U}(\rho) &= \int_1^\rho \rho\, \dot{U}(\rho)\, \frac{\textrm{d}\tau}{\tau} \\
& =\int_1^\rho \tau\, \dot{U}(\tau)\, \frac{\textrm{d}\tau}{\tau} \;=\; U(\rho)  - \,U(1) \;>\, 0
\end{split}
\end{equation*}
Now the following limit is easily seen to exist for $\,h \in \Ha (\A, \ast)\,$:
\begin{equation}
 \infty \,>\, \dot{U} (1)  =  \lim_{\rho \searrow 1} \;\rho\,\dot{U} (\rho) \; \geqslant 0\,,\;\;\;\;\textrm{in particular}\,,\;\;\; \;\dot{\mathscr R}(1) \geqslant 0\,.
\end{equation}
\begin{definition} The so defined quantity $\dot{\mathscr R}(1) = \dot{\mathscr R_h} (1)\,$ is called the \textit{initial speed}, it represents the initial data of \textit{harmonic evolution of circles} under the mapping $\,h\,$.
\end{definition}
We just found that the initial speed is always nonnegative for $h \in\mathcal H(\mathbb A,\ast)$; later in Proposition \ref{SchII}  we shall see that it equals 1 for conformal mappings.
The first derivative $\dot{U}(\rho)$, now  well defined in the left closed interval $[1, R)$,  is continuous. Finally, formula~\eqref{eq25} remains valid for $\varepsilon = 1 $ and reads as,
\begin{equation}\label{eq111}
\frac{1}{\pi} \iint_{A(1,\, \rho)} \abs{Dh}^2 \;=\;\rho\,\dot{U}(\rho) \; -\; \dot{U}(1)  \; <\;\infty\,,\;\;\textrm{for} \;\;1\leqslant \;\rho  < R
 \end{equation}
 This is what we wished to establish.\\

Other circular means require assuming that $\,h\,$ has continuous normal and tangential derivatives at the unit circle. This technical obstacle will be overcome by the following approximation argument.

\begin{lemma}\label{lemappr}\textnormal{\textsc{(weak approximation)}}  Given  $h\in \Ha (\A, \ast)\,,\,\mathbb A = A(1,R)\,,$  there exist annuli   $\mathbb A_k  = A(1,R_k) \,$  and harmonic mappings $h^k \in \Ha (\A_k, \ast)\,,$ $ \,k =1,2,...,$ such that
\begin{itemize}
\item {each $\,h^k\,$ is $\mathscr C^\infty-$smooth up to the inner boundary (the unit circle)}
\item {each $\,h^k\, : \mathbb T \overset{\textnormal{\tiny{onto}}}{\longrightarrow}\mathbb T$ is a homeomorphism homotopic to the identity}
\item {$\, 1<R_k \nearrow R\,$}
\item{$
\underset{k \rightarrow \infty}{\lim} h^k = \,h\,,\;\textrm{weakly in every Sobolev spaces}\; \mathscr W^{1,2}(\mathbb A_\circ\,\,\mathbb C),\,\\ {\;} \quad \quad \quad \quad \quad \;\;\textrm{where} \; \mathbb A_\circ = A(1, R_\circ)\;\textrm{and} \; 1 <R_\circ < R\,.$}\label{wlim}
\end{itemize}

\end{lemma}

\begin{remark}
Because of harmonicity the sequence $\{h^k\}$ and its all derivatives converge uniformly on every circle $\T_\rho=\{z \in \C \colon \abs{z}=\rho\}$, $1< \rho<R$.
\end{remark}
\begin{proof}
 Choose and fix a sequence of radii $\,r_1>r_2>...>1,$ close enough to 1 to have all circles $ \{ \xi ; \, |\xi| = r_k \}$ contained in $h(\mathbb A)\,.$ Since $h \colon \A \to \C$ is a real-analytic diffeomorphism, the level sets $\Gamma_k =  \{z \,; |h(z)|  =  r_k\,\}$
are real-analytic Jordan curves in $\mathbb A$. They separate the boundary components.
 Consider the doubly connected domain $\Delta_k \subset \A$ whose inner boundary is $\Gamma_k$ and outer boundary is that of $\,\mathbb A\,$. There is uniquely determined circular annulus $\,\mathbb A_k = A(1, R_k)\,$ which admits a conformal map $\Phi_k \colon \mathbb A_k     \overset{\textnormal{\tiny{onto}}}{\longrightarrow} \Delta_k \,.$    Such a map is $\mathscr C^\infty$ - smooth and injective  up to the closer of  $\,\mathbb A_k\,$. Actually, it extends conformally beyond the boundary of $\,\mathbb A_k\,$. Relevant details can be found in~\cite[p. 14]{Dubo} and~\cite[p. 41]{Pob}. We may assume that $\,\Phi_k\,$    preserves the order of boundary components, that is  $ \Phi_k (z) \in \Gamma _k\,,$ for $\,|z| = 1\,$. The outer radius of $\mathbb A_k$ is less than $\,R\,$. In fact, we have
 $$
   \log R_k  = \,\textrm{Mod}\, \Delta_k \,\; < \,\textrm{Mod} \,\mathbb A\; =\; \log R\,,\;\;\textrm{and} \;\;R_k \nearrow R
 $$
 Passing to a limit with a suitable subsequence of $\{\Phi_k\}$ results in a conformal map of $\,\mathbb A\,$ onto itself. Furthermore, composing $\,\Phi_k\,$ with a rotation of the $z$-variable (appropriately chosen for each function $\,\Phi_k\,$), we manage that
$\Phi_k(z)\,\rightarrow \, z\,$, for every $\,1 < \abs{z} <R $.  Actually the convergence is uniform together with all derivatives on every compact subset of the annulus $\,\mathbb A\,$.
 We  are now ready to define the harmonic mappings in question
\[h^k(z)= \frac{1}{r_k} h \big(\Phi_k (z)\big) \, , \qquad \mbox{ for } 1 < \abs{z} < R_k\]

To complete the proof of the lemma, we need only establish the existence of the weak limit in~\eqref{wlim}.
Consider an annulus $\,\mathbb A_\circ  = \{z\,;\, 1<|z|< R_\circ < R\}\,$ with outer boundary $\,\mathbb T_{R_\circ}\subset \mathbb A\,$. It is evident, by topology,  that $\Phi_k(\mathbb A_\circ)$ lies inside the Jordan curve $\Phi_k (\mathbb T_{R_{\circ}})\,$. Even more, there is an annulus $A(1,\,\rho)\,$, $\,1<\rho<R\,,$ such that $\Phi_k (\mathbb T_{R_\circ}) \subset A(1, \rho)\,$ for all $k=1,2,...$     . This is  because $\Phi_k$ converge uniformly on $\mathbb T_{R_{\circ}}$ to the identity map. Thus, again by topology,  $\Phi_k(\mathbb A_\circ) \subset A(1,\,\rho)$. Now, in view of~\eqref{eq111}, the Dirichlet energy of $h^k$ on $\mathbb A_\circ$ is free from $\,k\,$.  Indeed, by a conformal change of variables,  we see that

$$
\iint_{\mathbb A_\circ} \abs{Dh^k(z)}^2 \, \dtext z = r_k^{-2}\iint_{\Phi_k (\mathbb A_\circ )} \abs{Dh(z)}^2 \, \dtext z \, \leqslant \;
\iint_{ A(1,\, \rho) } \abs{Dh(z)}^2 \,\dtext z \; < \infty
$$

The proof of the lemma is completed by noting that $\{h^k\}$ converges weakly to $h$ in the Sobolev space $\mathscr W^{1,2}(\mathbb A_\circ, \mathbb C)$.
%These are later conveniently renumerated by $1,2,\dots$.
\end{proof}

\section{Some Related Results and More Preliminaries} %\label{Sec3}
Here is the first taste and sample of the utility of circular averages~\eqref{means}. We use them for conformal mappings to prove
Theorem~\ref{Schottky}. Although the result is classical and shorter proofs can be given (\cite[p.~333]{Neb} or ~\cite{BPC}),
this proof and the comments following it indicate in some detail the underlying strategy to be used for harmonic mappings.

\subsection{Schottky's theorem redeveloped} \label{Sec3.1}

\begin{proposition}\label{SchII}  The initial speed  of conformal evolution of circles always equals 1. That is, if a homeomorphism $\, h\in\mathscr H(\mathbb A, \ast)\,$ is conformal, then
\begin{equation*}
\dot{\mathscr R_h} (1) \;=\;\lim_{\rho\searrow 1} \frac{\textrm{d}}{\textrm{d}\rho}\left(\dashint_{\T_\rho} \abs{h}^2\right)^{\frac{1}{2}} \, = 1
\end{equation*}
Moreover, for each $\,1<\rho < R\,,$ the circular means of $\,|h|^2\,$ satisfy:
\begin{equation}\label{3ineq}
  U(\rho)= \dashint_{\T_\rho} \abs{h}^2 \,  \geqslant \rho^2,\;\;\;\;\;\dot{U}(\rho) \geqslant 2\rho\,,\;\;\;\;\;\; \ddot{U}(\rho) \geqslant\, 2.
 \end{equation}
  Equality occurs, somewhere at $\,\rho\in(1,R)$,  if and only if $\,h(z) = \lambda z\;, \;|\lambda| = 1\,$.
\end{proposition}

The classical Schottky theorem follows by imposing the  outer boundary condition, $\,|h(z)| = R_\ast \,$ for $\, |z| = R\,$, to infer that  $\, R_\ast  \geqslant R\,$. This can be reversed via consideration of the inverse conformal map, ascertaining Theorem \ref{Schottky}.
\begin{proof}  One compelling motive for studying the averages such as $\,U(\rho)\,$  is to take advantage of convexity properties of holomorphic (later harmonic) functions. Let us begin with the Laurent expansion
\begin{equation}
 h(z) \,= \, \sum_{n\in\mathbb Z}  a_n z^n \,, \;\;\;\;\;\;\;\;\; 1< |z| <  R\,.
\end{equation}
The system $\,\{z^n\}_{z\in\mathbb Z}\,$ is orthogonal on every circle $\,\mathbb T_\rho\,, \;  1<\rho < R\,$. Thus,
$$
 U(\rho) \; =\; \sum_{n\in\mathbb Z}  |a_n|^2\, \rho^{2n}\,, \;\;\;\;\;\;\;\;\; 1< \rho <  R\,.
$$

All that matters is to find a certain second order differential operator $\,\mathcal L\,: \mathscr C^2(1,R) \rightarrow \mathscr C(1,R)\,,$  acting on $\,U,$  that fits into the following
scenario:
 \begin{equation}\mathcal L[U] \geqslant 0\,,\; \textit{with equality if and only if} \;\;\;h(z) = \lambda z\,,\;\;|\lambda| = 1\, .
 \end{equation}

 Much of the essential properties of $\,U(\rho)\,$ are contained in such inequality when combined with the topological behavior  of $\,h\,$ near the inner boundary of  $\,\mathbb A\,$.
 The interested reader may wish to consult   ~\cite{IKO} for a fuller discussion of this topic.
 Here, for conformal case, direct computation shows that
 \begin{equation}\label{19}
  \mathcal L[U] \,:=\; \frac{1}{\rho} \,\frac{\textrm{d}}{\textrm{d}\rho}\left[   \rho^3\,\frac{\textrm{d}}{\textrm{d}\rho}\left(\frac{U}{\rho^2} \right ) \right] \; = \;4\sum_{n\in\mathbb Z} n(n-1) | a_n|^2 \rho^{2n-2} \geqslant  0
 \end{equation}
 Hence
 \begin{equation}
 \rho^3\,\frac{\textrm{d}}{\textrm{d}\rho}\left(\frac{U}{\rho^2} \right) \; \geqslant \;   \rho^3\,\frac{\textrm{d}}{\textrm{d}\rho}\left(\frac{U}{\rho^2} \right) _{|\,\rho =1}=\; \dot{U}(1) \,-\, 2\,U(1)\, =\; \dot{U}(1) \,-\,2 \,.
 \end{equation}
 Now the Cauchy-Riemann equations in polar coordinates
 \begin{equation}\label{CR}
     h_\rho \; =\; \frac{i}{\rho} \, h _\theta\,
 \end{equation}
 and topology (winding number along the inner boundary) come into play.  But first note that $\,h\,$ is $\,\mathscr C^1$-regular up to the inner boundary of $\mathbb A$. Even more, since $\,|h(z)| \equiv 1\,$  on $\, \mathbb T\,$, it extends as a conformal map slightly inside the unit circle.  Let us reveal in advance that at this point of the investigation of harmonic mappings we shall justify the assumption of  $\,\mathscr C^1$ - regularity  by  Lemma \ref{lemappr}.
 This makes it legitimate to perform the following computation,
 $$
 \dot{U}(1) \, =\; 2 \,\re \,\dashint_{\T} \bar{h}\,h_\rho \;=\; 2 \,\im \,\dashint_{\T} \bar{h}\,h_\theta \;=\; 2 \,\im \,\dashint_{\T} \frac{h_\theta}{h} = \; \frac{1}{\pi} \underset{|z|=1}{\Delta}\textrm{Arg} \,h(z) \;= 2\,.
 $$
 We just proved that  $\dot{\mathscr R}_h(1) = 1\,$. Returning to~\eqref{19} we infer that the function  $\, \rho\rightarrow \rho^{-2} \,U(\rho) \,$ is nondecreasing, and hence $\,U(\rho) \geqslant \rho^2\,$. We actually have slightly stronger inequality $\frac{\textrm{d}}{\textrm{d}\rho}\left(\frac{U}{\rho^2} \right) \; \geqslant \; 0\,$, which yields $\,\dot{U} \geqslant \frac{2}{\rho} \,U \geqslant  2\,\rho\,$.  Then it follows again from~\eqref{19} that $\,\ddot{U} \geqslant \frac{1}{\rho} \,\dot{U} \,\geqslant \, 2\,$.

 Finally, if equality occurs in one of~\eqref{3ineq} for some $\,1 < \rho < R\,$,  we infer from~\eqref{19}
that $\,a_n = 0 \,$, except for $\,a_0\,$ and $\, a_1\,$, which gives a linear function  $\,h(z) =  a_0 + a_1 z\, $. Since $\,|h(z)| \equiv 1\,$ on $\,\mathbb T\,$ we conclude that $\, h(z)= \lambda z\,,\;|\lambda | = 1\,$, as desired.
\end{proof}

Let us emphasize the principal features of this proof and indicate possible generalizations.

\subsection{Laplace equation} First of all, the Cauchy-Riemann system~\eqref{CR} was fundamental in the above proof. For, it let us replace the radial derivative $\,h_\rho\,$ by its angular derivative $\,\frac{i}{\rho} h_\theta\,$; whereas for harmonic mappings our strategy must rely on the second order equation:
\begin{equation}\label{Lapl}
h_{\rho \rho} \;=\; -\,\frac{\,h_{\rho}\,}{\rho} \;-\;\frac{h_{\,\theta \theta}\, }{\rho^2}
\end{equation}
\subsection{Orthogonal components}\label{Sec 3.3}
The basic complex harmonic functions in the annulus $\A$ are the integer powers $z^n$, $\bar z^n$, $n=0, \pm 1, \pm 2, \dots $, and the logarithm $\log \abs{z}$. If we combine these functions suitably in pairs, we obtain an analogue of the Laurent expansion for harmonic functions
\begin{equation}\label{hfeq1}
h(z)= \sum_{n\in \mathbb Z} h_n (z)\;,\;\;\;\;\;\;U(\rho) \,=\, \dashint_{\mathbb T_\rho}|h|^2= \sum_{n\in \mathbb Z} U_n(\rho) =\sum_{n\in \mathbb Z}\; \dashint_{\mathbb T_\rho}|h_n|^2
\end{equation}
where $h_n(z)= a_n z^n+b_n \bar z^{-n}$ for $n \ne 0$ and $h_0 (z)= a_0 \log \abs{z}+b_0$. These components are orthogonal on every circle $\mathbb T_\rho = \{z \colon \abs{z}=\rho\}$, $1 \le \rho < R$. We refer to $\,a_n\,$ and $\,b_n\,$ as the Fourier coefficients of $\,h\,$.

%Let us record for later use the following formulas:
%\begin{equation}\label{hfeq2}
 %U_n(\rho) = \begin{cases}  \;\;|a_n|^2 \rho^{2n}\;  + \;\;2\,\re \, (\bar{a}_n b_n)\;+ \;\;|b_n|^2 \rho^{-2n}\;, \;\;\;\textrm{for} & n  %\ne 0 \\
 %\;|a_0|^2 \log^2\rho \;+\,2\,\re \, (\bar{a}_0 b_0) \log\rho\, + \,|b_0|^2 \;, \;\;\;\;\textrm{for}& n = 0 .
%\end{cases}
%\end{equation}

\subsection{Nitsche mappings} The leading terms, corresponding to $\,n=1\,$, will take control of the geometric behavior of the mappings $\,h\in\mathcal H(\mathbb A,\ast)$. We conveniently scale them to become the identity map on the unit circle. So obtained functions, called the \textit{Nitsche maps}, form a one parameter family:
\[
\hbar_v(z)= \frac{1}{2}\left( z+ \frac{1}{\overline{z}}\right)\;+ \; \frac{v}{2}\left( z- \frac{1}{\overline{z}}\right)\; \;,\; \;1\leqslant |z|\,<\infty \, , \quad 0 \leqslant v < \infty.
\]
It is worth noting that, in the above sum, the terms have vanishing Neumann and Dirichlet  boundary values at the inner circle, respectively.
Moreover,  the parameter $\,v\,$ represents the initial speed of the evolution of circles,
$$
 \frac{\textrm{d}}{\textrm{d}\rho}\left(\dashint_{\T_\rho} \abs{\hbar_v}^2\right)^{\frac{1}{2}}_{|\,\rho =1} \, = v \,,\;\;\;\;\;\;\;\hbar_v\,\colon \,A(1,R) \rightarrow A(1,R_\ast)
$$
The outer radius of the image of any annulus $\, A(1,R)\,$ complies with the conjectural Nitsche bound:
\begin{equation}\label{outR}
 R_\ast \; = \frac{1}{2}\left( R+ \frac{1}{R}\right)\; +\;\,\frac{v}{2}\left( R- \frac{1}{R}\right)
\;\geqslant  \,\frac{1}{2}\left( R+ \frac{1}{R}\right)
\end{equation}
Of particular significance is the \textit{critical Nitsche map} with zero initial speed,
\[
\hbar(z) = \hbar_\circ(z)= \frac{1}{2}\left( z+ \frac{1}{\overline{z}}\right)\,.
\]

\subsection{Extreme properties of the Nitsche mappings}
It has been proven in ~\cite{AIM, IO} that $\,\hbar_v\,$ is exactly the mapping, unique up to rotation, that minimizes the Dirichlet energy subject to all homeomorphisms $\, h\in \mathscr H(\mathbb A,\mathbb A^*)\,$ between annuli $\mathbb A = A(1,R)\,$ and $\,\mathbb A^{\ast}\, =  A(1, R_{\ast})$, with $\,R_\ast\,$ given by formula~\eqref{outR}.

For another extreme property, we consider a subclass of homeomorphisms $\,h\in \mathcal H^\circ_v(\mathbb A,*) \subset \mathcal H(\mathbb A,*)\,$ (harmonic evolutions of circles) with  initial speed $\,\dot{\mathscr R}_h(1) \geqslant  v\,$
  and having zero first moment on the inner boundary; precisely,  $\,b_0 = 0\,$ in the expansion~\eqref{hfeq1}. Such is the mapping  $\,\hbar_v $. In this class of harmonic mappings the following generalization of the Nitsche bound has been established in ~\cite{IKO}:

\begin{equation*}
\mathscr R_h(R) \,=\; \left(\dashint_{\T_R} \abs{h}^2\right)^{\frac{1}{2}} \geqslant \frac{1}{2}\left( R+ \frac{1}{R}\right)\; +\;\,\frac{v}{2}\left( R- \frac{1}{R}\right) \;,\;\;\;h\in \mathcal H^\circ_v( \mathbb A,*)
\end{equation*}
Equality occurs for $\,h = \hbar_v\,$, uniquely up to a rotation. The idea of the proof is very much the same as that for Proposition \ref{SchII}. However, one has to remodel the operator $\mathcal L\,$, say $\,\mathcal L = \mathcal L_v\,$ in order to obtain $\,\mathcal L_v[U] = 0\,$ for $\,h=\hbar_v\,$.

The precondition $\,b_0=0\,$ seems to be redundant regardless of the initial speed. In fact this has already been removed for $\,\textrm{Mod}\,\mathbb A\,$ sufficiently small, depending on the initial speed $\,v\,$, by more advanced estimates in ~\cite{IKO}.
Unexpectedly, the redundancy of this seemingly minor precondition turned out to be the key difficulty in completing the proof of the Nitsche conjecture.

\subsection{The operator $\,\mathcal L(U)$} To each  Nitsche map $\,\,\hbar_v, \;0\leqslant v < \infty\,$, there corresponds a differential operator which tells us something about convexity properties of integral means of harmonic maps, with $\,\hbar_v\,$ as extreme case. We have already seen such an operator for $\,v=1\,$ in the proof of Proposition \ref{SchII}.  But our interest here is in the critical Nitsche map $\,\hbar (z) =  \frac{1}{2} ( z+ \overline{z}^{\,-1})\,$ of initial speed $\,v=0\,$. Associated with $\,\hbar\,$  is its operator $\,\mathcal L\,$, see  ~\cite{IKO}, which we introduce here in three different ways as follows:

\begin{align} %\label{equ3}
\mathcal L [U]\,&:= \;\ddot{U}+ \frac{3-\rho^2}{\rho (\rho^2+1)} \dot{U} - \frac{8}{(\rho^2+1)^2}U  \label{L1} \\
&=\frac{\rho^2+1}{\rho^3} \frac{\dtext}{\dtext \rho} \left[ \rho^3 \frac{\dtext}{\dtext \rho}  \left(\frac{U}{\rho^2+1}  \right)  \right] \label{L2}\\
&=\dashint_{\T_\rho} \left[  2 \abs{h_\rho}^2 + \frac{2}{\rho^2} \abs{h_\theta}^2 - 2 \frac{\rho^2-1}{\rho (\rho^2+1)} \abs{h^2}_\rho - \frac{8}{(\rho^2+1)^2} \abs{h}^2  \right]
\label{L3}
\end{align}

The interested reader may wish to learn the following principal feature of this operator; all the orthogonal components $\,h_n(z) = a_n z^n +  b_n z ^{-n}\,,\;$ with $\, n = \pm 1, \pm 2, ... , $ give rise to nonnegative values of $\,\mathcal L[U]\,$, whereas $\,\mathcal L[U]\equiv 0\,$ for $\,\hbar\,$. We shall not appeal to these properties explicitly, though they suggest us how to approach sharp estimates. But for the completeness of our arguments we still owe to the reader a verification of formula~\eqref{L3}; the one in~\eqref{L2} is clear. It is interesting to note that~\eqref{L3} contains no second derivatives of $\,h\,$.
Let us begin with the integral means and their derivatives
\begin{align*}
U(\rho) &= \dashint_{\T_\rho} \abs{h}^2 \\
\dot{U}(\rho) &= \dashint_{\T_\rho} \abs{h^2}_\rho = 2 \dashint_{\T_\rho} \re \bar h\,  h_\rho \\
\ddot{U} (\rho) &=   2 \dashint_{\T_\rho} \re \left( \bar h_\rho \, h_\rho \, + \, \bar h\,  h_{\rho \rho} \right)
\end{align*}
At this and only this  stage  we use the Laplace equation~\eqref{Lapl}, which yields:
\[\ddot{U} (\rho) = \dashint_{\T_\rho} \left[ 2 \abs{h_\rho}^2 - \frac{1}{\rho} \abs{h^2}_\rho   \right] - \frac{2}{\rho^2} \re \dashint_{\T_\rho} \bar h \, h_{\theta \theta} \]
The latter term, upon integration by parts along the circle $\T_\rho$, becomes $\dashint_{\T_\rho} \abs{h_\theta}^2$. In this way we represent $\ddot{U}(\rho)$ by using only  first derivatives of $h$,
\[\ddot{U}(\rho) = \dashint_{\T_\rho} \left[ 2 \abs{h_\rho}^2 - \frac{1}{\rho} \abs{h^2}_\rho + \frac{2}{\rho^2} \abs{h_\theta}^2  \right]\]
Substituting these formulas for $U$, $\dot{U}$ and $\ddot{U}$ into~\eqref{L1} we arrive at~\eqref{L3}.

\section{The Case $\Mod \A \le 1$}\label{thin}
The generalized Nitsche bound in~\eqref{genbound} is a straightforward consequence of a rather sophisticated identity for harmonic functions. We created this identity for the convenience of the reader in order to capture the essentials of the proof of  Theorem~\ref{th2} when $\,\textrm{Mod} \,\mathbb A \,\leqslant 1$.
\begin{proposition}\label{proplesse}
Let $h$ be a complex harmonic function in the annulus $\A=A(1,R)$, $1<R< \infty$, that is $\mathscr C^1$-smooth
up to the boundary. Then
\begin{equation}\label{theidentity}
\begin{split}
 & \hskip-0.2cm \frac{2R^2}{R^2+1} \dashint_{\T_R} \abs{h}^2  - \frac{R^2+1}{2} \dashint_{\T} \abs{h}^2\\
& -(R^2-1) \dashint_{\T} \abs{h} \abs{h}_\rho - (R^2-1) \log R\;  \dashint_{\T} \im \bar h \left(h_\theta - i h\right) \\
%& \hskip0.1cm =\\
&  = \frac{1}{\pi} \iint_{\A} \left[(R^2-1) \log \frac{R}{\rho} + \frac{R^2-\rho^2}{\rho^2}\right] \cdot \left| \frac{\rho h_\rho -i h_\theta}{1+\rho^2} - \frac{2 \rho^2\, h}{(1+\rho^2)^2}  \right|^2 \\
& \hskip0.4cm + \frac{1}{\pi} \iint_{\A} \left[(R^2-\rho^2) - (R^2-1) \log \frac{R}{\rho}  \right] \cdot \left| \frac{\rho h_\rho +i h_\theta}{1+\rho^2} + \frac{2 \, h}{(1+\rho^2)^2}   \right|^2
\end{split}
\end{equation}
\end{proposition}
The derivation of this identity is postponed until the end of this section.

\subsection{Proof of Theorem~\ref{th2} when $\Mod \A \le 1$}\label{subsec55}$\;$\\

No restriction for $R$ is needed in Proposition \ref{proplesse}. Nevertheless, for the proof of the Nitsche bound we shall have to restrict the outer radius to the interval  $1< R \le e$ in order to ensure that  double integrals of~\eqref{theidentity} are nonnegative. The first integrand is certainly positive for every $1 \le \rho \le R$. However, the second integral needs the restriction $1 \le \rho \le R \le e$. We have
\begin{equation*}
(R^2-\rho^2)-(R^2-1) \log \frac{R}{\rho} \ge 0 \, , \quad \mbox{whenever } 1 \le \rho \le R \le e
\end{equation*}
To see this, note that the expression in the left hand side represents a concave function in $\,\rho$ with nonnegative values at the endpoints of the interval [1, R].

For the proof of~\eqref{genbound} we  take $\sigma \in (1,R)$, then choose and fix a radius $R_0$ between $\,\sigma\,$ and $\,R\,$,  $\,1<\sigma<R_0<R$. Consider the annulus $\A_0=  A(1,R_0) \subset \A$. Given $h\in \Ha (\A, \ast)$ we appeal to Lemma~\ref{lemappr} to construct harmonic homeomorphisms $h^k\in \Ha (\A_0, \ast)$ that are $\mathscr C^1$-smooth up to the inner boundary of $\A_0$  and converge to $h$ weakly in $\mathscr W^{1,2}(\A_0)$.  Before proceeding to the  identity~\eqref{theidentity} we note three particulars concerning integral means over the unit circle:
\begin{itemize}
\item[(i)]{$\displaystyle \dashint_{\T} \abs{h^k}^2=1$, \;\; because $\abs{h^k} \equiv 1$ on $\T$}\\
\item[(ii)]{$\displaystyle  \dashint_{\T} \abs{h^k} \abs{h^k}_\rho \ge 0$ }
\end{itemize}
This is because homeomorphisms $h^k$ take circles $\T_\rho$, $1<\rho<R_0$ into Jordan curves inside which there lies the unit disk. Precisely, we have at every point of the unit circle:
\[\abs{h^k}_\rho = \lim_{\rho \searrow 1} \frac{\abs{h^k (\rho e^{i \theta})}-1 }{\rho -1} \ge 0\]
We also have the identity
\begin{itemize}
\item[(iii)] $ \displaystyle  \im \dashint_{\T} \overline{h}^k \left( h^k_\theta - i h^k\right)=0$
\end{itemize}
which follows from the computation of the winding number of $h^k$ around $\T$,
\[\dashint_{\T} \overline{h}^k h^k_\theta = \dashint_{\T} \frac{h^k_\theta}{h^k}=i= i \dashint_{\T} \abs{h^k}^2\]
On substituting (i-iii) into~\eqref{theidentity}, with $\,R\,$ replaced by $\,\sigma\,$,  we obtain
\[
\begin{split}
 & \hskip-0.2cm \frac{2\sigma^2}{\sigma^2+1} \dashint_{\T_{\sigma}} \abs{h^k}^2  - \frac{\sigma^2+1}{2} \\
&  \ge \frac{1}{\pi} \iint_{A(1, \sigma)} \left[(\sigma^2-1) \log \frac{\sigma}{\rho} + \frac{\sigma^2-\rho^2}{\rho^2}\right] \cdot \left| \frac{\rho h^k_\rho -i h^k_\theta}{1+\rho^2} - \frac{2 \rho^2\, h^k}{(1+\rho^2)^2}  \right|^2 \\
& \hskip0.4cm + \frac{1}{\pi} \iint_{A(1, \sigma)} \left[(\sigma^2-\rho^2) - (\sigma^2-1) \log \frac{\sigma}{\rho}  \right] \cdot \left| \frac{\rho h^k_\rho +i h^k_\theta}{1+\rho^2} + \frac{2 \, h^k}{(1+\rho^2)^2}  \right|^2
\end{split}
\]
We are going to pass to the limit as $k \to \infty$. Note that $h^k\rightrightarrows h$ uniformly on $\T_{\sigma}$ and weakly in  $\mathscr W^{1,2}(\A_0)$. Passing to the limit in the double integrals results in the desirable estimate from below, due to  lower semicontinuity of the double integrals,  in which  $h^k, h^k_\rho$ and $h^k_\theta$ converge weakly in $\mathscr L^2(\A_0)$.
\begin{equation}\label{589436737856}
\begin{split}
 & \hskip-0.2cm \frac{2\sigma^2}{\sigma^2+1} \dashint_{\T_{\sigma}} \abs{h}^2  - \frac{\sigma^2+1}{2} \\
&  \ge \frac{1}{\pi} \iint\limits_{A(1, \sigma)} \left[(\sigma^2-1) \log \frac{\sigma}{\rho} + \frac{\sigma^2-\rho^2}{\rho^2}\right] \cdot \left|  \frac{\rho h_\rho -i h_\theta}{1+\rho^2} - \frac{2 \rho^2\, h}{(1+\rho^2)^2} \right|^2 \\
& \hskip0.4cm + \frac{1}{\pi} \iint\limits_{A(1, \sigma)} \left[(\sigma^2-\rho^2) - (\sigma^2-1) \log \frac{\sigma}{\rho}  \right] \cdot \left| \frac{\rho h_\rho +i h_\theta}{1+\rho^2} + \frac{2 \, h}{(1+\rho^2)^2}  \right|^2
\end{split}
\end{equation}
Hence
\[ \frac{2\sigma^2}{\sigma^2+1} \dashint_{\T_{\sigma}} \abs{h}^2  - \frac{\sigma^2+1}{2} \ge 0 \]
or, equivalently
\begin{equation}\label{1437895}
\left(  \dashint_{\T_{\sigma}} \abs{h}^2 \right)^\frac{1}{2} \ge \frac{1}{2} \left(\sigma + \frac{1}{\sigma}\right)
\end{equation}
This is the estimate we wished to obtain.

For uniqueness statement in Theorem~\ref{th2}, we assume that equality occurs in~\eqref{1437895}. Then we see from~\eqref{589436737856} that $h$ must satisfy the following equations:
\begin{equation*}
\begin{cases}
%\begin{split}
& \left( \rho \frac{\partial}{ \partial \rho} - i \frac{\partial}{ \partial \theta} - \frac{2 \rho^2}{1+\rho^2}  \right) h =0 \\
& \left( \rho \frac{\partial}{ \partial \rho} + i \frac{\partial}{ \partial \theta} + \frac{2 }{1+\rho^2}  \right) h =0
%\end{split}
\end{cases}
\end{equation*}
Adding and subtracting the equations we uncouple the $\rho$ and $\theta$ derivatives,
\[
\begin{cases}
%\begin{split}
\rho \, h_\rho &= \frac{\rho^2-1}{\rho^2+1}h  \\
i\,  h_\theta &=h
%\end{split}
\end{cases}
\]
The general solution takes the form $h (\rho e^{i \theta}) = a \left(\rho + \frac{1}{\rho}\right) e^{i \theta}$ where $a$ is any complex number. Since $\abs{h}\equiv 1$ on $\T$ we conclude that $\abs{a}=1/2$.  The proof of Theorem~\ref{th2} in case $1< R \le e$ will therefore be accomplished  once we establish the identity~\eqref{theidentity}.\qed

\subsection{Proof of identity ~\eqref{theidentity}}
The identity is obtained by integrating $\,\mathcal L[U]\,$ against a weight over the interval $\,\{\rho : 1\, < \rho\,<R\}\,$. We shall integrate by parts using the divergence form of $\,\mathcal L[U]\,$ in~\eqref{L2}. The weight must be carefully crafted in order not to produce the derivative $\,\dot{U}(R)\,$ as a boundary term. This term is out of control; it can even be infinity in some cases, exactly when the energy of $\,h\,$ over $\,\mathbb A\,$ is infinite.  There is essentially only one weight that suits well to this conception, namely $\frac{\rho (R^2-\rho^2)}{\rho^2+1}$.
Since $\,\mathcal L[U]\,$ vanishes for the critical Nitsche mapping, it is natural to try to simplify computation by making a substitution
\[h(z)= \frac{1}{2} \left( z + \frac{1}{\bar z} \right) g(z)\, , \qquad 1 < \abs{z}<R  \]

In order to express $\LL[U]$ by means of $g$ we compute the terms under the integral sign in~\eqref{L3}
\begin{equation*}\begin{split}
\abs{h}^2&=\frac{(\rho^2+1)^2}{4\rho^2}\abs{g}^2 \\
\abs{h_{\theta}}^2&=\frac{(\rho^2+1}{4\rho^2}\abs{g_{\theta}+ig}^2 \\
\frac{d}{d\rho}\left(\frac{\rho h}{\rho^2+1}\right)&=\frac{e^{i\theta}}{2}g_{\rho}
\end{split}\end{equation*}
Therefore,
\begin{equation*}\begin{split}
\LL[U]&=\frac{(\rho^2+1)^2}{2\; \rho^4}\,
\dashint_{C_{\rho}}\left(\abs{g_{\theta}+ig}^2-\abs{g}^2+\rho^2\abs{g_{\rho}}^2\right) \\
&=\frac{(\rho^2+1)^2}{2\; \rho^2}\
\dashint_{C_{\rho}}\left[\;\abs{g_{\rho}}^2+\rho^{-2}\abs{g_{\theta}}^2+2\rho^{-2}\im (\bar g g_{\theta})\;\right]
\end{split}\end{equation*}
In this way we arrive  at somewhat simpler formula
\begin{equation*}%\label{givename}
\mathcal L [U] = \frac{(\rho^2+1)^2}{\rho^2} \dashint_{\T_\rho}  \left[  \abs{g_z}^2 +  \abs{g_{\bar z}}^2 + \frac{1}{\rho^2} \im \left(\bar g \, g_\theta \right) \right]
\end{equation*}
We shall now exploit the divergence form of $\mathcal L [U]$ in~\eqref{L2}.   Multiply both sides of~\eqref{L2} by the weight $\frac{\rho (R^2-\rho^2)}{\rho^2+1}$ and integrate  from $\rho=1$ to $\rho=R$.
\[
\begin{split}
\int_1^R\frac{R^2-\rho^2}{\rho^2}  \frac{\dtext}{\dtext \rho}  & \left[ \rho^3  \frac{\dtext}{\dtext \rho}  \left(  \frac{U}{\rho^2 +1} \right)  \right] = \int_1^R \;\frac{\rho\,(R^2-\rho^2)}{\,\rho^2+1\;}\;\mathcal L[U]  \;\;\;= \\ &\\&= \int_1^R \frac{(R^2-\rho^2)(\rho^2+1)}{\rho} \dashint_{\T_\rho}  \left[  \abs{g_z}^2 +  \abs{g_{\bar z}}^2 + \frac{1}{\rho^2} \im \left(\bar g \, g_\theta \right) \right]
\end{split}
\]
Integration by parts of the lefthand side leaves only boundary terms
\begin{equation}\label{equ4}
\begin{split}
\frac{2R^2}{R^2+1} U(R)&- \frac{R^2+1}{2}U(1)-\frac{R^2-1}{2}\dot{U}(1)\\
&= \int_1^R \frac{(R^2-\rho^2)(\rho^2+1)}{\rho} \dashint_{\T_\rho} \left(  \abs{g_z}^2 + \abs{g_{\bar z}}^2 \right)\\ &\; \; +  \int_1^R \frac{(R^2-\rho^2)(\rho^2+1)}{\rho^3} \dashint_{\T_\rho}\im \left(  \bar g \, g_\theta \right)
\end{split}
\end{equation}
Here we split the righthand side for the purpose of integrating only the second term by parts. To achieve this objective, we represent the factor in front of $\dashint_{\T_\rho} \im \bar g \, g_\theta\;$ as  derivative of a function that vanishes at $\rho=R$, namely:
\[\frac{(R^2-\rho^2)(\rho^2+1)}{\rho^3}= \frac{\dtext}{\dtext \rho} \left[  \frac{(R^2-\rho^2)(\rho^2-1)}{2\rho^2} - \left(R^2-1\right) \log \frac{R}{\rho}   \right]\]
Integration by parts produces additional boundary term and formula~\eqref{equ4} becomes
\begin{equation}\label{equ5}
\begin{split}
\frac{2R^2}{R^2+1} & \dashint_{\T_R} \abs{h}^2  - \frac{R^2+1}{2} \dashint_{\T}\abs{h}^2 - \frac{R^2-1}{2} \dashint_{\T}\abs{h^2}_\rho \\ & =  \int_1^R \frac{(R^2-\rho^2)(\rho^2+1)}{\rho} \dashint_{\T_\rho}\left( \abs{g_z}^2 +  \abs{g_{\bar z}}^2\right) \\
& - \int_1^R \left[  \frac{(R^2-\rho^2)(\rho^2-1)}{2\rho^2} - \left(R^2-1\right) \log \frac{R}{\rho}    \right] \frac{\dtext}{\dtext \rho} \dashint_{\T_\rho} \im \left( \bar g \, g_\theta \right)\\
&+ (R^2-1) \log R \, \dashint_{\T} \im \left( \bar g \, g_\theta \right)\
\end{split}
\end{equation}

With the aid of commutation rule ~\eqref{comrule} we pass the $\rho$-differentiation inside the mean integral and then integrate by parts along the circle $\,\mathbb T_\rho\,$ to obtain

\[
\begin{split}
\frac{\dtext}{\dtext \rho} \dashint_{\T_\rho} \im (\bar g \, g_\theta) &= \im \dashint_{\T_\rho} \left( \bar g_\rho g_\theta + \bar g g_{\rho \theta} \right) \\
&= \im \dashint_{\T_\rho} \left( \bar g_\rho g_\theta - \bar g_\theta  g_{\rho} \right) = 2 \rho \dashint_{\T_\rho} \left(  \abs{g_z}^2 - \abs{g_{\bar z}}^2\right)
\end{split}
\]
Now, the righthand side of~\eqref{equ5} takes the form
\[ \begin{split}
&\int_1^R  \left[\frac{(R^2-\rho^2)(\rho^2+1)}{\rho} - \frac{(R^2-\rho^2)(\rho^2-1)}{\rho} + (R^2-1) 2 \rho \log \frac{R}{\rho}  \right] \dashint_{\T_\rho} \abs{g_z}^2\\
&  + \int_1^R  \left[\frac{(R^2-\rho^2)(\rho^2+1)}{\rho} + \frac{(R^2-\rho^2)(\rho^2-1)}{\rho} - (R^2-1) 2 \rho \log \frac{R}{\rho}  \right] \dashint_{\T_\rho} \abs{g_{\bar z}}^2\\
&\hskip2cm  +(R^2-1) \log R \dashint_{\T} \im \left( \bar g\,  g_\theta \right)
\end{split}
\]
which simplifies to
\[
\begin{split}
\frac{1}{\pi} \int_1^R & \left[ \frac{R^2-\rho^2}{\rho^2} + (R^2-1) \log \frac{R}{\rho} \right] \int_{\T_\rho} \abs{g_z}^2
\\ &+  \frac{1}{\pi} \int_1^R  \left[ R^2-\rho^2 - (R^2-1) \log \frac{R}{\rho} \right] \int_{\T_\rho} \abs{g_{\bar z}}^2\\& \hskip2cm+ (R^2-1) \log R \;\dashint_{\T} \im \left( \bar g\,  g_\theta \right)
\end{split}\]
Finally, the entire formula~\eqref{equ5} reads as
\begin{equation*}
\begin{split}
\frac{2R^2}{R^2+1}  \dashint_{\T_R} \abs{h}^2 &  - \frac{R^2+1}{2} \dashint_{\T}\abs{h}^2 - (R^2-1) \dashint_{\T}\abs{h}\abs{h}_\rho \\ & -   (R^2-1) \log R\; \dashint_{\T} \im \left( \bar g\,  g_\theta \right)\\
&=  \frac{1}{\pi} \iint_{\A}  \left[ \frac{R^2-\rho^2}{\rho^2} + (R^2-1) \log \frac{R}{\rho} \right]   \abs{g_z}^2 \\
& +  \frac{1}{\pi} \iint_{\A}  \left[ (R^2-\rho^2) - (R^2-1) \log \frac{R}{\rho} \right]   \abs{g_{\bar z}}^2
\end{split}
\end{equation*}
To conclude with the identity~\eqref{theidentity} it only remains to observe that
\[
\begin{split}
\im \left( \bar g \, g_\theta \right) &= \frac{4 \rho^2}{(1+\rho^2)^2} \im \big(\bar h h_\theta-i\bar hh\big) \\
\abs{g_z} &=  \left|\frac{\rho h_\rho-ih_\theta}{1+\rho^2} - \frac{2 \rho^2h}{(1+\rho^2)^2}  \right|  \\
\abs{g_{\bar z}} &=  \left| \frac{\rho h_\rho+ih_\theta}{1+\rho^2} - \frac{2 h}{(1+\rho^2)^2}  \right|
\end{split}
\]
This ends the proof of Proposition~\ref{proplesse}.  \qed

\begin{remark}
Recall from Section~\ref{Sec 3.3} the orthogonal decomposition $h(z)= \sum_{n\in \mathbb Z} h_n (z)\;$ of a harmonic function, where $h_n(z)= a_n z^n+b_n \bar z^{-n}$ for $n \ne 0$ and $h_0 (z)= a_0 \log \abs{z}+b_0$. Putting this in the identity~\eqref{proplesse} breaks it up into independent identities, one for each term $\,h_n$. This could be another way of proving Proposition~\ref{proplesse}. Precisely, for each individual term $\,h_n$ the problem reduces to showing that so obtained quadratic form with respect to the complex variables $\,a_n, b_n\in \mathbb C\,$ is nonnegative. However, it would not make the proof simpler; one immediately encounters the difficulty in the simplest possible case $\,h(z)\,\equiv 1\,$. Try it! Nevertheless, the orthogonal decomposition and the method of quadratic forms for the Fourier coefficients of $\,h\,$   will prove extremely useful for the case $\Mod \A \ge 1$,
see Proposition \ref{bigrprop} in which a counterpart of the identity~\eqref{proplesse} cannot be even formulated without appeal to the coefficients $\,a_n, b_n\,$.
\end{remark}

\section{The Case $\Mod \A \ge 1$} %\label{thick}

Let us return to the initial conditions ~\eqref{in1}--\eqref{in3} discussed in the introduction. An example  below demonstrates that, in contrast to the case $\Mod \A\le 1$, condition~\eqref{in3} cannot be omitted when $\Mod \A$ is large. This is the underlying reason why the method of Section~\ref{thin} could not work for  arbitrary values of $\Mod \A$.
\subsection{The example}
\begin{example}\label{ex}
Fix $0<a<1$ and let $\lambda$ be a positive number to be chosen later. Define for $\,|z| > 1\,$,
\[
h(z)=\frac{1+a\bar z}{\bar z+a}+\lambda\log\abs{z},\quad \textrm{thus}\;\;\; \abs{h(z)}\le 1+\lambda\log\abs{z}
\]
It is clear that $h$ is harmonic in $\{z\colon \abs{z}>1\}$ and satisfies~\eqref{in1} with $r=1$.
We represent $\,h\,$ as a sum of two terms which are orthogonal on every circle $\T_\sigma$, $\sigma\ge 1$,
\[
h(z)=(a+\lambda\log\abs{z})+\frac{1-a^2}{\bar z+a}\;,
\]
In this way we accomplish  the following computation
\[
\frac{\dtext}{\dtext\sigma}\dashint_{\T_{\sigma}}\abs{h}^2 =
\frac{\dtext}{\dtext\sigma}(a+\lambda\log \sigma)^2
+\frac{\dtext}{\dtext\sigma}\dashint_{\T_{\sigma}}\frac{(1-a^2)^2}{\abs{\bar z+a}^2}
\]
It follows that~\eqref{in2} holds if $\lambda$ is chosen to be sufficiently large, depending only on the parameter $a$, which remains fixed. However, the generalized Nitsche bound~\eqref{genbound} fails on circles of large radius, because $\abs{h(z)}$  exhibits logarithmic growth as $\,|z| \rightarrow \infty\,$.
The reader may wish to note that the average Jacobian determinant over the unit circle is independent of $\,\lambda\,$ and is negative,  $\,\dashint_{\T} \det Dh  = -\,\dashint_{\T} |z+a|^{-4} |\textrm{d}z| = - \frac{1+ a^2}{(1-\,a^2)^3}$.
\end{example}

In view of the above example, we need to take advantage of the additional information that the averaged Jacobian is nonnegative at the unit circle, as stated in  condition ~\eqref{in3}.
\subsection{An inequality for harmonic functions}

The condition ~\eqref{in3} will come into play via the following inequality.

\begin{proposition}\label{bigrprop}
Let $h\colon A(1,R)\to\C$ be a harmonic function that is $\mathscr C^1$-smooth up to the inner boundary circle $\T$. Denote by $f\colon \overline{\DD}\to\C$ the harmonic extension of $h$ to the closed unit disk. Then for all $\sqrt{7}\le \rho< R$ we have
\begin{equation}\label{trueforall}
\begin{split}
\dashint_{\T_\rho}\abs{h}^2 & -\left(\frac{\rho+\rho^{-1}}{2} \right)^2\dashint_{\T}\im (\bar h h_\theta) \\
& - 2\dashint_{\T} \abs{h}\abs{h}_\rho   -\frac{\rho^2-4-\rho^{-2}}{2}\dashint_{\T} J_h\\ &  -\frac{\rho^2-4-\rho^{-2}}{4\pi } \left[\int_{\T} \det D\!f-  \iint_{\mathbb D} \abs{D\!f}^2\right]\ge 0
\end{split}
\end{equation}
\end{proposition}

Since $\sqrt{7}<e$, Proposition~\ref{bigrprop} covers all values $e \le \rho < R$.\\
\proof
Our proof of~\eqref{trueforall} involves the Fourier coefficients of $h$. These are complex numbers
$a_n,b_n$, $n\in\Z$, that appear in the orthogonal expansion
\[
h(z)=\sum_{n\in\Z}h_n(z)=a_0\log\abs{z}+b_0+\sum_{n\ne 0}(a_n z^n + b_n \bar z^{-n})
\]

First observe that the harmonic extension of $h$ inside the unit disk is expressed by two infinite sums
\begin{equation*}%\label{fseries}
f(z)= \sum_{n \ge 0} (a_n +b_n)z^n +  \sum_{n < 0} (a_n +b_n)\bar z^n.
\end{equation*}
This is certainly true for mappings $h\in \Ha (\A,\ast)$ that are continuous  up to the inner boundary. Then
the terms in~\eqref{trueforall} can be computed using orthogonality of the powers of $z=\rho\, e^{i \theta}$.
\begin{equation}\label{mess}
\begin{split}
&\dashint_{\T_\rho} \abs{h}^2 = \abs{a_0 \log \rho +b_0}^2 + \sum_{n \ne 0} \abs{a_n\rho^n + b_n \rho^{-n}}^2;\\
&\dashint_{\T} \abs{h}\abs{h}_\rho = \frac{1}{2} \dashint_{\T} \abs{h^2}_\rho =  \re( a_0\bar b_0)+ \sum_{n \ne 0} n(\abs{a_n}^2-\abs{b_n}^2);
\\
&\dashint_{\T} \im \left( \bar h \, h_\theta \right) = \sum_{n \ne 0} n \abs{a_n+b_n}^2; \hskip0.5cm     \dashint_{\T} J_h = \sum_{n \ne 0} n^2(\abs{a_n}^2-\abs{b_n}^2) ;\\
&\dashint_{\T} \det D\!f =  \sum_{n \ne 0} n\abs{n} \abs{a_n+b_n}^2; \hskip0.5cm
\iint_{\mathbb D} \abs{D\!f}^2 \;=\, 2\,\pi\sum_{n \ne 0} \abs{n} \abs{a_n+b_n}^2
\end{split}
\end{equation}
 Upon substituting these terms to~\eqref{trueforall} the lefthand side becomes a quadratic form $\, Q = \,Q\binom{\,\dots\,,\, a_{-1},\, \,a_0,\, a_1,\,\dots}{\,\dots\,,\, b_{-1},\,\,\, b_0,\,\, b_1,\,\dots}\,$ with respect to the complex variables $\,a_n,\,b_n\,\in\mathbb C\,, n \in \mathbb Z\,$,
which we aim to show to be nonnegative
for every complex numbers $a_n,b_n$, $n\in\Z$. More precisely, $Q$ splits into an infinite sum of quadratic forms, each of which depends only
on two complex variables,
\begin{equation}\label{qform1}
Q=\sum_{n\in\Z} Q_n(a_n,b_n),\qquad a_n,b_n\in\C,
\end{equation}
where
\begin{equation*}%\label{qform2}
Q_n(\xi,\zeta)=A_n(\rho)\abs{\xi}^2+B_n(\rho)\abs{\zeta}^2\,+\,2\,C_n(\rho)\re (\xi \bar \zeta).
\end{equation*}
For example,
\begin{equation}%\label{qform3}
Q_0(\xi,\zeta)=\abs{\xi\log\rho+\zeta}^2-2\re(\xi\bar\zeta)
\end{equation}
which is positive definite as long as $\,\log\rho > \frac{1}{2}\,$.
A key to the extremal case lies in the property of the quadratic form $\,Q_1\,$, which is positive semidefinite,
\begin{equation*}%\label{qform4}
Q_1(\xi,\zeta)=\frac{(\rho^2-1)^2}{4\rho^2}\abs{\xi-\zeta}^2 \,> 0 \,, \;\;\textrm{unless}\;  \xi = \zeta\;,
\end{equation*}
 meaning that $\,Q_1(a_1,b_1) = 0\;$ yields $\,  h_1(z) =  a \left(z + \frac{1}{\bar z} \right)\,$.
 The other quadratic forms will be shown to be positive definite.
\begin{lemma}\label{techlem}
For every $n \ne 0,1$ and $\xi, \zeta \in \C$ we have
\begin{equation*}%\label{qform5}
Q_n(\xi,\zeta)>0\;, \qquad \text{unless } \ \xi=\zeta=0\, .
\end{equation*}
\end{lemma}
We postpone the proof of this technical lemma until Appendix 1.

\subsection{Jacobian--Energy Inequality: Proof of Theorem~\ref{th3}}
The second inequality in~\eqref{th3form} is immediate from
\[\abs{D\!f(z) }^2 = 2\left( \abs{f_z}^2 + \abs{f_{\bar z}}^2 \right) \ge  2\left( \abs{f_z}^2 - \abs{f_{\bar z}}^2 \right)  = 2\,\det D\!f(z)\, \]
Integration of the Jacobian determinant over the unit disk gives the area of the  image, which equals~$\pi$.

To prove the first part of~\eqref{th3form}, it suffices to consider the case when $f$ is sense-preserving.
At every point $z=e^{i\theta}$ of the unit circle $\T$ we have $\abs{f(z)}=1$
and $\overline {f(z)} f_\theta(z)=i\abs{f_\theta(z)}$. Now we compute the Jacobian determinant at $z\in \T$ as follows.
\begin{equation*}%\label{jfid}
\det D\!f(z) =\im (\overline {f_\rho} f_\theta)
=\im \big[(\overline{ f_\rho} f) (\overline{f} f_\theta)\big]=\abs{f_\theta}\re (\overline{ f_\rho} f)=\abs{f_\theta}\cdot\abs{f}_{\rho}\;.
\end{equation*}
Since $f$ is harmonic, we have $\Delta (\abs{f}^2)=2\abs{D\!f}^2$. Green's formula yields
\begin{equation*}%\label{enerf}
\iint_\DD \abs{D\!f}^2 = \frac{1}{2}\int_{\T} \abs{f^2}_\rho = \int_{\T} \abs{f}_\rho.
\end{equation*}
Thus our goal is to prove
\begin{equation}\label{jacthm1}
\int_{\T} \abs{f}_\rho (\abs{f_\theta}-1) \,d\theta\;\ge 0
\end{equation}
Let us write $f(e^{i\theta})=e^{i\xi(\theta)}$, where $\xi \in \mathscr C^1([0, 2\pi])$ is increasing and satisfies $\xi(2\pi)= \xi(0)+ 2\pi$.
By Poisson's formula, for $0\le \rho< 1$ and $0\le \theta\le 2\pi$, we have
\begin{equation*}%\label{poi1}
f(\rho e^{i\theta})=\frac{1}{2\pi}\int_0^{2\pi}\frac{1-\rho^2}{1-2\rho\cos(\theta-\phi)+\rho^2}\,e^{i\xi(\phi)}\,d\phi
\end{equation*}
Multiply by $e^{-i\xi(\theta)}$ and take the real part:
\begin{equation*}%\label{poi2}
\re\left(\overline{f(e^{i\theta})}f(\rho e^{i\theta})\right)
=\frac{1}{2\pi}\int_0^{2\pi}\frac{1-\rho^2}{1-2\rho\cos(\theta-\phi)+\rho^2}\,\cos\big[\xi(\theta)-\xi(\phi)\big]\,d\phi.
\end{equation*}
Now the normal derivative $\,\abs{f}_{\rho}\,$ at the unite circle can be computed as follows.
\begin{equation*}%\label{poi3}
\lim_{\rho\nearrow 1} \frac{1-\re(\overline{f(e^{i\theta})}f(\rho e^{i\theta}))}{1-\rho}
=\frac{1}{2\pi}\int_0^{2\pi}\frac{1-\cos\big[\xi(\theta)-\xi(\phi)\big]}{1-\cos(\theta-\phi)}\,d\phi.
\end{equation*}
The lefthand side of~\eqref{jacthm1} becomes a double integral with respect to both $\theta$ and $\phi$
\begin{equation*}%\label{poi4}
2\pi \int_{\T} \abs{f}_{\rho}(\abs{f_\theta}-1)\;d\theta \,
= \int_0^{2\pi}\!\!\int_0^{2\pi} \frac{1-\cos\big[\xi(\theta)-\xi(\phi)\big]}{1-\cos(\theta-\phi)}\,\left(\xi'(\theta)-1\right)\;d\theta\;d\phi\,
\end{equation*}
We are reduced to showing that
\begin{lemma} \label{Le5}
Suppose $\xi \in \mathscr C^1([0, 2\pi])$ is increasing and $\,2\pi\,$ periodic; that is  $\,\xi(2\pi) = \xi(0) + 2\pi\,$. Then
\begin{equation}\label{poi55}
\int_0^{2\pi}\!\!\int_0^{2\pi} \frac{1-\cos\big[\xi(\theta)-\xi(\phi)\big]}{1-\cos(\theta-\phi)}\,\left(\xi'(\theta)-1\right)\;d\theta\;d\phi\, \;\geqslant 0
\end{equation}
 Equality occurs if and only if  $\,\xi(\theta) = \theta + \const\,$.
\end{lemma}
 Elementary, though lengthy proof of this lemma is given in Appendix 2.

\subsection{Proof of Theorem~\ref{th2} when $\,\Mod \A \ge 1$} Since the result of Section~\ref{thin} applies to the restriction of $h$ to $A(1,e)$, we only concern ourselves with integral means of $h$ over the circles of radius $\rho \ge e$.   First assume that $h\in \Ha(\A,\ast)$ is $\mathscr{C}^1$-smooth up to the inner boundary. Thus we may compute the winding number of $h$ around $\T$,
\[
\dashint_{\T} \im \bar h h_\theta = \im \dashint_{\T}\frac{h_\theta}{h} = 1,
\]
 Moreover,
\[
\dashint_{\T} \abs{h}\abs{h}_\rho \ge 0,\quad \text{because } \ \abs{h}_\rho\ge 0 \ \text{ on } \ \T
\]
and
\[
\dashint_{\T} \det Dh \ge 0, \quad \text{because } \ \det Dh \ge 0 \ , \;\text{pointwise}
\]
Let $f \colon \overline \DD \rightarrow \C $ be the continuous extension of $h_{|_{\T}}$ that is harmonic in $\DD$. Since $h(z)-f(1/\bar z)$ is harmonic in $A(1,R)$ and vanishes on $\T$, it has an extension as a harmonic function on  $A(1/R,R)$. In particular, $f \in \mathscr C^1(\overline{\DD})$.  By Theorem~\ref{th3} we have
\[
\int_{|z|=1} \det D\!f (z) \; \ge \iint_{|z|\leqslant 1 } \abs{D\!f (z)}^2
\]
Substituting these inequalities  to  ~\eqref{trueforall} yields the desired Nitsche bound
\begin{equation*}%\label{desired}
\left(\dashint_{\T_\rho}\abs{h}^2\right)^{1/2}\ge \frac{1}{2}\left(\rho+\frac{1}{\rho}\right).
\end{equation*}
There is no difficulty to relax the $\mathscr C^1$-smoothness assumption.  Let us go into this in detail with the aid of Lemma~\ref{lemappr}, in much the same way as in Section~\ref{subsec55}.  We  fix $\sigma \in (1,R)$, then choose and fix a radius $R_0$ so that $1<\sigma<R_0<R$. Consider the annulus $\A_0=  A(1,R_0) \subset \A$. Given $h\in \Ha (\A, \ast)$ we appeal to Lemma~\ref{lemappr} to construct harmonic homeomorphisms $h^k\in \Ha (\A_0, \ast)$ that are $\mathscr C^1$-smooth up to the boundary of $\A_0$  and converge to $h$ weakly in $\mathscr W^{1,2}(\A_0)$.
 This time we have four particulars concerning integral means over the unit circle:
\begin{itemize}
\item[(i)]{$\displaystyle \;\dashint_{\T}\; \abs{h^k}^2=1$ }\\
\item[(ii)]{$\displaystyle  \;\dashint_{\T}\; \abs{h^k}\,  \abs{h^k}_\rho \ge 0$ }\\
\item[(iii)]{$\displaystyle \; \dashint_{\T}\; \im \bar h^k h^k_\theta  = 1$ }\\
\item[(iv)]{$\displaystyle  \; \dashint_{\T}\; \det Dh^k\ge 0$ }
\end{itemize}
 For each  $h^k$ the Nitsche bound holds,
\begin{equation}\label{7757587}
\left(\dashint_{\T_\rho} \abs{h^k}^2  \right)^\frac{1}{2} \ge \frac{1}{2} \left(\rho + \frac{1}{\rho}\right)\, , \quad k=1,2,\dots
\end{equation}
It is essential that before passing to the limit we have ignored, based on (i)-(iv),  the integral means of $h$ and their derivatives over the inner boundary of $\,\mathbb A\,$.
Now, because of harmonicity the sequence $\{h^k\}$ converges uniformly on every circle $\T_\rho$, $1< \rho <R\,,$ but not necessarily for $\,\rho = 1\,$.
Passing to the limit we obtain the Nitsche bound for every $h \in \Ha (\A, \ast)$. \\
The uniqueness statement is somewhat delicate, it takes into account a stronger variant of~\eqref{7757587}, which in terms of the quadratic forms in~\eqref{qform1} reads as
\begin{equation*}%\label{later}
\dashint_{\T_\rho} \abs{h^k}^2- \frac{1}{4} \left(\rho + \frac{1}{\rho}\right)^2 \ge \sum_{|n| \leqslant N} Q_n(a_n^k,b_n^k) \,,\;\;\;\;\;\textrm{for every}\;N = 1,2, ...
\end{equation*}
Here $a_n^k,\, b_n^k$ are the associated Fourier coefficients of $h^k$; the quadratic forms with $\,|n| > N\,$,  being positive definite,  are omitted.
 We can now pass to the limit
\begin{equation}\label{1664}
\dashint_{\T_\rho} \abs{h}^2 - \frac{1}{4}  \left(\rho  + \frac{1}{\rho}\right)^2  \ge  \sum_{|n|\leqslant N} Q_n (a_n, b_n) \ge 0
\end{equation}
Passage to the limit in the finite sum of quadratic forms  is justified because for every fixed $n$ we have $\,\lim_{k\to\infty} a_n^k = a_n$ and $\,\lim_{k\to\infty} b ^k_n = b_n$.

Finally, for the uniqueness, suppose that $\;\dashint_{\T_\rho} \abs{h}^2 = \frac{1}{4} \left(\rho + \frac{1}{\rho}\right)^2$,  where $\,1< \rho <R$. Then~\eqref{1664} yields  $Q_n(a_n,b_n)= 0\,$ for each integer $n$. Hence, by Lemma \ref{techlem}, all the coefficients
$a_n$ and $b_n$ vanish except for the case $n=1$. This leaves us with the only possible function
\[h(z)=az+ \frac{a}{\bar z}= \frac{1}{2} \left(z+ \frac{1}{\bar z}\right) e^{i\alpha}\]
for some $0 \le \alpha < 2\pi$, because $\abs{h(z)}=1$ on $\T$. The proof of Theorem~\ref{th2} is complete, modulo Lemma~\ref{techlem}. \qed

\section{Appendix 1, Proof of Lemma~\ref{techlem}} %\label{secforms}
The quadratic form  $Q_n(\xi,\zeta)=A_n(\rho)\abs{\xi}^2+B_n(\rho)\abs{\zeta}^2\,+\,2\,C_n(\rho)\re (\xi \bar \zeta)$  in the decomposition~\eqref{qform1} is obtained by putting $\, h(z) = h_n(z) = a_n z  + b_n \bar{z} ^{-1}\,$  into the left hand side of~\ref{trueforall}.
We first do the case of positive indices; that is,\\

{\bf Case  $n \ge 2$.}  With the aid of formulas~\eqref{mess} we find the coefficients:
\begin{equation*}%\label{coef1}
\begin{split}
A_n&=\rho^{2n}-\frac{n}{4}(\rho+\rho^{-1})^2-2n-\frac{2n^2-n}{2}(\rho^2-4-\rho^{-2}); \\
B_n&=\rho^{-2n}-\frac{n}{4}(\rho+\rho^{-1})^2+2n+\frac{n}{2}(\rho^2-4-\rho^{-2});  \\
C_n&=1-\frac{n}{4}(\rho+\rho^{-1})^2-\frac{n^2-n}{2}(\rho^2-4-\rho^{-2}). \\
\end{split}
\end{equation*}
We need to show that $A_n$ and $B_n$ are positive and $A_nB_n >C_n^2$.
Ignoring the term $\rho^{-2n}$ in $B_n$, we obtain the estimate
\begin{equation}\label{Bnest1}
\begin{split}
B_n&\ge \frac{n\rho^{2}}{4}\left(-(1+\rho^{-2})^2+8\rho^{-2}+2(1-4\rho^{-2}-\rho^{-4})\right)
\\ &=\frac{n\rho^{2}}{4} (1 -2\rho^{-2}-3\rho^{-4}) \ge \frac{n\rho^{2}}{4}\left(1 -\frac{2}{7}-\frac{3}{49}\right)\ge \frac{n\rho^2}{7}
\end{split}
\end{equation}
Next, estimate $A_n$ from below  as follows.
\begin{equation}\label{Anest1}
\begin{split}
A_n&=\rho^{2n}-\left(n^2-\frac{n}{4}\right)\rho^2+\left(4n^2-\frac{9n}{2}\right)
+\left(n^2-\frac{3n}{4}\right)\rho^{-2} \\
&\ge \rho^{2n}-\left(n^2-\frac{n}{4}\right)\rho^2
\end{split}
\end{equation}
Regarding $C_n$, note that $C_n\le 0$ for all $n \ge 2$, and
\begin{equation}\label{Cnest1}
\begin{split}
\abs{C_n}&=-C_n=\left(\frac{n^2}{2}-\frac{n}{4}\right)\rho^2+\left(\frac{5n}{2}-2n^2-1\right)+\left(\frac{3n}{4}-\frac{n^2}{2}\right)\rho^{-2} \\
&\le \left(\frac{n^2}{2}-\frac{n}{4}\right)\rho^2+\left(\frac{5n}{2}-2n^2-1\right)
\end{split}
\end{equation}

With $n=2$, inequality~\eqref{Anest1} yields $A_2\ge \rho^4(1-7\rho^{-2}/2)\ge \rho^4/2$, which together with~\eqref{Bnest1}
and~\eqref{Cnest1} imply
\[
\begin{split}
A_2B_2-C_2^2 &\ge \frac{1}{7}\rho^6 -\left(\frac{3\rho^2}{2}-4\right)^2 =
\left(\frac{1}{7}\rho^4-\frac{9}{4}\rho^2+12-16\rho^{-2}\right)\rho^2\\ &> \left(\frac{1}{7}\rho^4-\frac{9}{4}\rho^2+9\right)\rho^2>0
\end{split}
\]
The latter inequality  holds for all $\rho>0$ because $\frac{36}{7}>\frac{81}{16}$.
Thus $Q_2$ is positive definite.

When $n\ge 3$, we simplify~\eqref{Cnest1}  further by ignoring a negative term  as follows
\begin{equation}\label{Cnest2}
\begin{split}
\abs{C_n}& \le \left(\frac{n^2}{2}-\frac{n}{4}\right)\rho^2+\left(\frac{5n}{2}-2n^2-1\right)
\le \left(\frac{n^2}{2}-\frac{n}{4}\right)\rho^2 \\
&=\left(\frac{1}{2}-\frac{1}{4n}\right)n^2\rho^2
\end{split}
\end{equation}
In light of~\eqref{Bnest1} and~\eqref{Cnest2} the inequality $A_nB_n > C_n^2$ will follow once we show that
\begin{equation}\label{goal1}
A_n > 7 \left( \frac{1}{2} - \frac{1}{4n}  \right)^2 n^3 \rho^2
\end{equation}
To this end, we use~\eqref{Anest1}
\begin{equation}\label{Anest2}
\begin{split}
A_n& \ge \rho^{2n}-\left(n^2-\frac{n}{4}\right)\rho^2
\ge \left(\frac{\rho^{2n-2}}{n^3}-\frac{1}{n}\right) n^3\rho^2\\
& \ge \left(\frac{7^{n-1}}{n^3}-\frac{1}{n}\right)n^3\rho^2
\end{split}
\end{equation}
When $n=3$, a direct computation shows that~\eqref{Anest2} implies~\eqref{goal1}. When $n\ge4$, we use the fact that $7^{n-1}n^{-3}$ is increasing in $n$ to obtain
\[ A_n \ge \left(\frac{7^3}{4^3}-\frac{1}{4}\right)n^3\rho^2>\left(2-\frac{1}{4}\right)n^3\rho^2
=\frac{7}{4}n^3\rho^2,\]
from which~\eqref{goal1} follows.

{\bf Case  $n\le -1$.} For convenience we set  $n=-m$, where $m$ is a positive integer.  With this new notation the coefficients of $Q_{n}$ are
\begin{equation}\label{coef4}
\begin{split}
A_{-m}&=\rho^{-2m}+\frac{m}{4}(\rho+\rho^{-1})^2+2m+\frac{m}{2}(\rho^2-4-\rho^{-2}); \\
B_{-m}&=\rho^{2m}+\frac{m}{4}(\rho+\rho^{-1})^2-2m+\frac{2m^2+m}{2}(\rho^2-4-\rho^{-2}); \\
C_{-m}&=1+\frac{m}{4}(\rho+\rho^{-1})^2+\frac{m^2+m}{2}(\rho^2-4-\rho^{-2})
\end{split}
\end{equation}
Organizing in powers of $\rho$, we find
\begin{equation}\label{coef5}
A_{-m} = \rho^{-2m}+ \frac{3m}{4}\rho^2 +\frac{m}{2} -\frac{m}{4}\rho^{-2} \ge \frac{3m}{4}\rho^2
\end{equation}
Similarly,
\begin{equation}\label{coef5b}
\begin{split}
B_{-m}&=\rho^{2m}+\left(m^2+\frac{3m}{4}\right)\rho^2-\left(4m^2+\frac{7m}{2}\right)-\left(m^2+\frac{m}{4}\right)\rho^{-2} \\
&\ge\rho^{2m}+\left(m^2+\frac{3m}{4}\right)\rho^2-\left(5m^2+4m\right)
\end{split}
\end{equation}
It is clear by~\eqref{coef4} that $C_{-m}$ is positive, so
\begin{equation}\label{coef5c}
\begin{split}
\abs{ C_{-m}} &=\left(\frac{m^2}{2}+\frac{3m}{4}\right)\rho^2 +\left(1-\frac{3m}{2}-2m^2\right)
-\left(\frac{m^2}{2}+\frac{m}{4}\right)\rho^{-2} \\
&\le \left(\frac{m^2}{2}+\frac{3m}{4}\right)\rho^2+\left(1-\frac{3m}{2}-2m^2\right)
\end{split}
\end{equation}
With $m=1$, inequalities~\eqref{coef5}--\eqref{coef5c} yield
\[
\begin{split}
A_{-1}B_{-1}-C_{-1}^2 &\ge \frac{3\rho^2}{4}\left(\frac{11\rho^2-36}{4}\right)-\left(\frac{5\rho^2-10}{4}\right)^2\\ &
=\frac{8\rho^2(\rho^2-1)-100}{16}\ge \frac{8\cdot 7\cdot 6-100}{16}>0
\end{split}
\]

When $m\ge 2$, we ignore the last term in~\eqref{coef5c} and obtain
\begin{equation}\label{coef6c}
\abs{C_{-m}} \le \left(\frac{m^2}{2}+\frac{3m}{4}\right)\rho^2\le \frac{7}{8}m^2\rho^2
\end{equation}
In light of~\eqref{coef5} and~\eqref{coef6c} the inequality $A_{-m}B_{-m} > C^2_{-m}$ will follow once we  show that
\begin{equation}\label{coef6}
B_{-m}>\frac{49}{48}m^3\rho^2,\qquad m=2,3,\dots
\end{equation}
For this we return to~\eqref{coef5b}. Since $\rho^2\ge 7$, it follows that
\begin{equation*}%\label{coef7}
\begin{split}
\frac{B_{-m}}{m^3\rho^2} &\ge  \frac{\rho^{2m-2}}{m^3}+\left(\frac{1}{m}+\frac{3}{4m^2} \right)
-\left(\frac{5}{m}+\frac{4}{m^2} \right)\rho^{-2} \\
&\ge \frac{7^{m-1}}{m^3}+\frac{2}{7m}+\frac{5}{28m^2}
\end{split}
\end{equation*}
where the latter is the minimum value of the former expression in~$\rho$, attained at $\rho^2=7$. It equals $\frac{17}{16}$ when $m=2$, while for $m\ge 3$ we have
\[\frac{B_{-m}}{m^3\rho^2} \ge \frac{7^{m-1}}{m^3}\ge \frac{49}{27}>\frac{49}{48},\]
proving~\eqref{coef6}. This completes the proof of Lemma~\ref{techlem}.\qed

\section{Appendix 2, Proof of Lemma~\ref{Le5}} %\label{Lem5}
The double integral in Lemma \ref{Le5} is the set $\,\mathbb Q = [0, 2\pi] \times [0, 2\pi]\,$.
We can write $\xi(\theta)=\theta+\zeta$, where
$\zeta=\zeta(\theta)$ is a $\mathscr C^1$-smooth $2\pi$-periodic function. Observe that
\begin{equation}\label{phidif}
\phi-\theta\le \zeta(\theta)-\zeta(\phi)\le 2\pi+\phi-\theta,\qquad 0\le \phi\le \theta\le 2\pi
\end{equation}
In particular, $\zeta'(\theta)\ge -1$ for all $\theta$. The left hand side of~\eqref{poi55} is equal to
\begin{equation}\label{poi5}
\iint_{\mathbb Q} \frac{1-\cos\big[\theta-\phi+\zeta(\theta)-\zeta(\phi)\big]}{1-\cos(\theta-\phi)}\,\zeta'(\theta)\;d\theta\, d\phi
\end{equation}
With the help of the trigonometric  identity $\cos(x+y)=\cos x\cos y-\sin x\sin y$ we represent~\eqref{poi5} as the sum $A+B$, where
\begin{equation*}%\label{int1}
\begin{split}
A&=\iint_{\mathbb Q} \frac{1-\cos(\theta-\phi)\cos\big[\zeta(\theta)-\zeta(\phi)\big]}{1-\cos(\theta-\phi)}\zeta'(\theta)\;d\theta\, d\phi\\
B&=\iint_{\mathbb Q} \frac{\sin(\theta-\phi)\sin\big[\zeta(\theta)-\zeta(\phi)\big]}{1-\cos(\theta-\phi)}\zeta'(\theta)\;d\theta\, d\phi.
\end{split}
\end{equation*}
Since $\zeta$ is periodic, the integral of $\zeta'$ over $\mathbb Q $ is equal to $0$. Let us subtract this integral from $A$ to write it as:
\begin{equation*}%\label{int2}
A=\iint_{\mathbb Q} \frac{\cos(\theta-\phi)\big[1-\cos(\zeta(\theta)-\zeta(\phi))\big]}{1-\cos(\theta-\phi)}\zeta'(\theta)\;d\theta\, d\phi.
\end{equation*}
We divide the set $\mathbb Q$ into
\[ \mathbb Q^+:=\{(\theta,\phi)\colon \cos(\theta-\phi)\ge 0\} \; \mbox{ and } \;  \mathbb Q^-:=\{(\theta,\phi)\colon \cos(\theta-\phi)< 0\}\] Accordingly, the integral $A$ splits into $A^+$ and $A^-$. The
$A^+$ part can be estimated from below using the inequality $\zeta'(\theta) \ge -1$,
\begin{equation}\label{int3}
\begin{split}
A^+&=\iint_{\mathbb Q^+} \frac{\cos(\theta-\phi)\big[1-\cos(\zeta(\theta)-\zeta(\phi))\big]}{1-\cos(\theta-\phi)}\zeta'(\theta)\;d\theta\, d\phi \\
&\ge \iint_{\mathbb Q^+} \frac{-\cos(\theta-\phi)\big[1-\cos(\zeta(\theta)-\zeta(\phi))\big]}{1-\cos(\theta-\phi)}\;d\theta\, d\phi
\end{split}
\end{equation}
In $A^-$ part we perform integration by parts with respect to $\theta$. No boundary terms appear because $\zeta$ is periodic and $\cos(\theta-\phi)$ vanishes on the common boundary of $\mathbb Q^+$ and $\mathbb Q^-$.
\begin{equation*}%\label{int4}
\begin{split}
A^-&=\iint_{\mathbb Q^-} \frac{\cos(\theta-\phi)(1-\cos\big[\zeta(\theta)-\zeta(\phi)\big]}{1-\cos(\theta-\phi)}\zeta'(\theta)\;d\theta\, d\phi \\
&=\iint_{\mathbb  Q^-} \frac{\cos(\theta-\phi)}{1-\cos(\theta-\phi)}\frac{d}{d\theta}\left\{
\zeta(\theta)-\zeta(\phi)-\sin\big[\zeta(\theta)-\zeta(\phi)\big]\right\}\;d\theta\, d\phi \\
&=\iint_{\mathbb Q^-} \frac{\sin(\theta-\phi)}{(1-\cos(\theta-\phi))^2}\left\{\zeta(\theta)-\zeta(\phi)-\sin\big[\zeta(\theta)-\zeta(\phi)\big]\right\}\;d\theta\, d\phi .
\end{split}
\end{equation*}
We also integrate $B$ by parts with respect to $\theta$
\begin{equation}\label{int5}
\begin{split}
B&=\iint_{\mathbb Q} \frac{\sin(\theta-\phi)}{1-\cos(\theta-\phi)}\frac{d}{d\theta}\left\{1-\cos\big[\zeta(\theta)-\zeta(\phi)\big]\right\} \;d\theta\, d\phi \\
&=\iint_{\mathbb Q} \frac{1-\cos\big[\zeta(\theta)-\zeta(\phi)\big]}{1-\cos(\theta-\phi)}\;d\theta\, d\phi
\end{split}
\end{equation}
The latter integral splits as $B^+ + B^-$, where $B^+$ is taken over $\mathbb Q^+$ and $B^-$ over $\mathbb Q^-$.
The sum $A^++B^+$ is estimated by combining~\eqref{int3} and~\eqref{int5}.
\begin{equation}\label{int6}
A^++B^+\ge \iint_{\mathbb Q^+} \left\{1-\cos\big[\zeta(\theta)-\zeta(\phi)\big]\right\}\,d\theta\,d\phi\ge 0
\end{equation}
Finally, we write the sum $A^-+B^-$ using shorthand notation $\alpha=\theta-\phi$ and $\beta=\zeta(\theta)-\zeta(\phi)$.
\begin{equation}\label{int7}
A^-+B^-= \iint_{\mathbb Q^-} \frac{(1-\cos \alpha)(1-\cos  \beta) +(\beta-\sin \beta)\sin \alpha}{(1-\cos \alpha)^2}\;d\theta\, d\phi
\end{equation}
The definition of $\mathbb Q^-$ and inequalities~\eqref{phidif} imply
\begin{equation}\label{abconstr}
\pi/2\le \alpha\le 3\pi/2\quad\text{ and } \quad -\alpha\le \beta\le 2\pi-\alpha
\end{equation}
We claim that the integrand in~\eqref{int7} is nonnegative, that is,
\begin{equation}\label{int8}
\Psi(\alpha,\beta):=(1-\cos \alpha)(1-\cos  \beta) +(\beta-\sin \beta)\sin \alpha \ge 0
\end{equation}
under the conditions~\eqref{abconstr}. From this~\eqref{poi55} will follow by adding up~\eqref{int6} and~\eqref{int7}.

It remains to verify ~\eqref{int8}.
Replacing the pair $(\alpha,\beta)$ with $(2\pi-\alpha,-\beta)$ if necessary, we may assume that $\pi/2\le \alpha\le \pi$.
Now, if $\beta\ge 0$, then $\beta-\sin\beta\ge 0$, and~\eqref{int8} follows. Suppose $\beta<0$.
Then $\Psi$ is increasing with respect to $\alpha \in [\pi/2,\pi]$. For a fixed $\beta\in [-\pi,0]$, the minimal admissible value of $\alpha$
under~\eqref{abconstr} is $\max(\pi/2,-\beta)$. This leads us to consider two cases.

\textit{\textbf{Case 1}}. If $-\pi/2\le \beta\le 0$, then
\begin{equation}\label{int9}
\Psi(\alpha,\beta)\ge \Psi(\pi/2,\beta)= 1-\cos  \beta + \beta-\sin \beta
\end{equation}
Differentiating the righthand side of~\eqref{int9}, we find that it is decreasing for $-\pi/2\le \beta\le 0$. Since it vanishes at $\beta=0$,
it follows that $\Psi(\alpha,\beta)\ge 0$.

\textit{\textbf{Case 2}}. If $-\pi\le \beta\le -\pi/2$, then
\begin{equation}\label{int10}
\begin{split}
\Psi(\alpha,\beta)&\ge \Psi(-\beta,\beta)=(1-\cos  \beta)^2 +\sin^2\beta -\beta\sin \beta \\
&=2-2\cos\beta -\beta\sin \beta
\end{split}
\end{equation}
Again, we find that the righthand side of~\eqref{int10} is decreasing for $-\pi\le \beta\le -\pi/2$. Its value at $\beta=-\pi/2$
is $2-\pi/2>0$. Thus $\Psi(\alpha,\beta)> 0$ in this case.

This completes the proof of~\eqref{poi55}. If equality holds in~\eqref{poi55}, then it also holds in~\eqref{int6}. The latter is
only possible if $\zeta$ is a constant function, i.e., when $\,\xi(\theta) = \theta + \const\,$. \qed

\bibliographystyle{amsplain}

\end{document}